\documentclass[aos,preprint]{imsart}

\RequirePackage[OT1]{fontenc}
\RequirePackage{amsthm,amsmath,amssymb}
\usepackage{epsfig}
\usepackage{graphicx}
\usepackage{epstopdf,color}
\usepackage[usenames,dvipsnames,svgnames,table]{xcolor}
\usepackage{multirow}
\usepackage{listings}
\RequirePackage{natbib}
\RequirePackage[colorlinks,citecolor=blue,urlcolor=blue]{hyperref}

\usepackage{xr}
\startlocaldefs
\numberwithin{equation}{section}
\theoremstyle{plain}
\newtheorem{lemma}{Lemma}[section]
\newtheorem{theorem}{Theorem}[section]

\newtheorem{remark}{Remark}[section]
\newtheorem{example}{Example}[section]
\endlocaldefs

\newcommand{\bu}{{\mathbf u}}

\newcommand{\x}{{\mathbf x}}

\newcommand{\z}{{\mathbf z}}

\newcommand{\CC}{\mathbb{C}}

\newcommand{\um}{\underline{m}}
\newcommand{\de}{\delta}

\newcommand{\tL}{\tilde L}
\newcommand{\R}{\tilde R}
\newcommand{\tr}{{\text{\rm tr}}}

\begin{document}

\begin{frontmatter}
\title{ High-dimensional covariance matrices in elliptical distributions with application to spherical test
}
\runtitle{ covariance matrices in elliptical distributions}

\begin{aug}
\author{\fnms{Jiang} \snm{Hu}\thanksref{t1,m1}
\ead[label=e1]{huj156@nenu.edu.cn}},

\author{\fnms{Weiming} \snm{Li}\thanksref{t2,m2}
\ead[label=e2]{li.weiming@shufe.edu.cn}},

\author{\fnms{Zhi} \snm{Liu}\thanksref{t3,m3}
\ead[label=e3]{liuzhi@umac.mo}}
\and
\author{\fnms{Wang} \snm{Zhou}\thanksref{t4,m4}
	\ead[label=e4]{stazw@nus.edu.sg}}

\thankstext{t1}{Jiang Hu's research is  partially supported by NSFC (No. 11771073) and the Fundamental Research Funds for the Central Universities.}
\thankstext{t2}{Weiming Li's research is partially supported by NSFC (No.
	11401037), MOE Project of Humanities and Social Sciences, No.\ 17YJC790057, and Program of IRTSHUFE.}
\thankstext{t3}{Zhi Liu's research
	is supported by FDCT127/2016/A3  of Macau and NSFC (No.11401607).}
\thankstext{t4}{Wang Zhou's research is partially supported by the MOE Tier 2 grant MOE2015-T2-2-039 (R-155-000-171-112)  at the National University of Singapore.}
\runauthor{J. Hu et al.}

\affiliation{Key Laboratory for Applied Statistics of MOE $\&$ Northeast Normal University\thanksmark{m1}, Shanghai University of Finance and Economics\thanksmark{m2}, University of Macau $\&$ UM Zhuhai Research Institute\thanksmark{m3}  and National University of Singapore\thanksmark{m4}}

\address{KLASMOE and School of Mathematics $\&$ Statistics\\ Northeast Normal University\\
	Changchun China 130024\\
\printead{e1}}

\address{School of Statistics and Management\\
	Shanghai University of Finance and Economics\\
	Shanghai China 200433\\
	\printead{e2}}

\address{Faculty of Science and Technology\\
	University of Macau\\
	Macau China\\
	\printead{e3}}

\address{Department of Statistics and Applied Probability\\ National University of Singapore\\
	Singapore 117546\\
	\printead{e4}}

\end{aug}

\begin{abstract}
This paper discusses fluctuations of linear spectral statistics of high-dimensional sample covariance matrices  when the underlying population follows an elliptical distribution.  Such population often possesses high order correlations among their coordinates, which have great impact on the asymptotic behaviors of linear spectral statistics.  Taking such kind of dependency into consideration, we establish a new central limit theorem for the linear spectral statistics in this paper for a class of elliptical populations. This general theoretical result has wide applications and, as an example, it is then applied to test the sphericity of elliptical populations.

\end{abstract}

\begin{keyword}[class=MSC]
\kwd[Primary ]{62H15}
\kwd[; secondary ]{ 62H10}
\end{keyword}

\begin{keyword}
\kwd{Covariance matrix}
\kwd{ High-dimensional data}
\kwd{ Elliptical distribution}
\kwd{ Sphericity test}.
\end{keyword}

\end{frontmatter}

\section{Introduction}
Large-scale statistical inference develops rapidly in the last two decades. This type of inference often relies on spectral statistics of certain random matrices in  high-dimensional frameworks, where both the dimension $p$ of the observations and the sample size $n$ tend to infinity.
Recall that   a {\em linear spectral statistic} (LSS) \citep{BSbook} of a $p\times p$ Hermitian random matrix $R_n$ is  of the form
\begin{equation}\label{lss}
\frac{1}{p}\sum_{i=1}^pf(\lambda_i)=\int f(x)dF^{R_n}(x),
\end{equation}
where $\lambda_1,\ldots,\lambda_p$ are the $p$ eigenvalues of $R_n$, $f$ is a function defined on $\mathbb R$, and $F^{R_n}=(1/p)\sum_{i=1}^p\delta_{\lambda_i}$ is called the {\em spectral distribution} (SD) of $R_n$.
Here $\delta_{a}$ denotes  the Dirac measure at the point $a$.
In \cite{LW02} and \cite{Schott07}, most test statistics are actually LSSs of sample covariance matrices. \cite{BJYZ09} made systematic corrections to several classical likelihood ratio tests to overcome the effect of high-dimension using LSSs of sample covariance matrices and F-matrices. Later, \cite{BJPZ15} derived the CLT for LSSs of a high-dimensional Beta matrix, which can be broadly used in multivariate statistical analysis, such as testing the equality of several covariance matrices, multivariate analysis of variance, and canonical correlation analysis, see   \cite{Anderson03}. Most recently, based on an LSS of regularized canonical correlation matrices, \cite{YP15} proposed a test for the independence between two large random vectors.   \cite{Gao16} applied LSSs  of sample correlation matrices
to the complete independence test for $p$ random variables and the equivalence test for $p$ factor loadings in a factor model.
Clearly, it is of great interests to investigate the behaviors of LSSs under various circumstances.

Specifically, let $\x_1,\ldots,\x_n$ be $n$ observations of $\x\in\mathbb R^p$, whose mean is zero and  covariance matrix is $\Sigma$. The sample covariance matrix is
$$
B_n=\frac{1}{n}\sum_{j=1}^n\x_j\x_j'.
$$
Our attention in this paper is focused  on the asymptotic properties of LSSs of $B_n$.
The earliest study on this problem dates back to  \cite{Jonsson82},  who obtained the
{\em central limit theorem} (CLT) for LSSs of  $B_n$ by assuming the population to be standard multivariate normal.  A remarkable breakthrough was done in \cite{BS04}, where the population is allowed to be a linear transform of a vector of independent and identically distributed (i.i.d.) random variables, i.e.,
\begin{equation}\label{linear}
\x=A\z.
\end{equation}
Here $A \in \mathbb R^{p\times p}$  is a non-random transformation matrix with $rank(A)=p$, and $\z=(z_1,\ldots,z_p)'$ with i.i.d.\ $z_i$'s satisfying
\begin{equation}\label{m-condition}
E(z_1)=0,\quad E(z_1^2)=1\quad\text{and}\quad E(z_1^4)=3.
\end{equation}
The fourth moment condition was later extended by \cite{PZ08} to $E(z_1^4)<\infty.$
Though these assumptions are fairly weak, their requirement of linearly dependent structure in \eqref{linear} still excludes a lot of important distributions. In particular, it excludes almost all distributions from the elliptical family \citep{FZ90}.

Elliptical distributions were originally introduced by \cite{K70} to generalize the multivariate normal distributions.  A random vector $\x$ with zero mean follows an elliptical distribution if and only if it has a stochastic representation \citep{FZ90}:
\begin{equation}\label{eds}
\x =  \xi A\bu,
\end{equation}
where the matrix $A \in \mathbb R^{p\times p}$ is non-random with $rank(A)=p$, $\xi\geq 0$ is a scalar variable representing the radius of $\x$, and  $\bu \in \mathbb R^{p}$ is the random direction, which is independent of $\xi$ and uniformly distributed on the unit sphere $S^{p-1}$ in $\mathbb R^p$,  denoted by $\bu\sim U(S^{p-1})$ in the sequel.
This family of distributions has been widely applied in many areas, such as statistics, economics and finance, which can describe fat (or light) tails and tail dependence among components of a population, see \cite{FZ90} and \cite{G13}. Evidently such distributions with high order correlations can not be modeled by the  linear transform model in \eqref{linear}.

A question raised immediately is that  how the nonlinear dependency affects the asymptotic behaviors of LSSs in high-dimensional frameworks?
Indeed, \cite{BZ08}  proved that the SD $F^{B_n}$ of $B_n$ converges to a common generalized Mar\v{c}enko-Pastur  law almost surely if,
for any sequence of symmetric matrices
$\{C_p\}$   with bounded spectral norm,
\begin{equation}\label{BZ-condition}
Var(\x' C_p\x)=o(p^2)
\end{equation}
as $p,n\to \infty$. This condition is also sharp for the convergence, see \cite{LY16} for an example. What is more,  this condition holds for  a list of well known elliptical distributions, such as multivariate normal distributions,  multivariate Pearson type II distributions,
power exponential distributions, and a more general family of multivariate Kotz-type distributions \citep{K75}, see Section \ref{sec:main} for more details. However, the limit of SD is not enough for  many procedures of statistical inference, such as confidence interval and hypothesis testing.
Therefore, in this paper, we will explore the fluctuations of LSSs of $B_n$, when the population belongs to elliptical distributions that satisfy the condition \eqref{BZ-condition}.
Compared with the pioneer work of \cite{BS04}, the main
difficulty of the current study lies in the fact that both the radius $\xi$ and direction $\bu$ introduce nonlinear dependence to the coordinates of the population $\x$, which can not be handled through the same way as they did for the linearly dependent structure.   Technically, we are facing the following three challenges. First, for paying the  cost of dropping linearly dependent structure, we have to add more moment conditions on  $\xi$, because  the finite fourth moment of $\xi/\sqrt{p}$ is no longer sufficient for the nonlinear dependence case (see  \eqref{m-con}). This is totally different from the linearly dependent structure case.  Second, we need to  figure out how the dependence of the entries of  $\xi\bu$ influences  the fluctuations of LSSs of $B_n$ (see  Remark \ref{rem2.3}). Third, we have to extend many fundamental conclusions in the independent case (\cite{BS04}) to accommodate  our non-linearly dependent structure; see Lemma  \ref{double-e}-\ref{lambda-bound} for example.  


The structure of this paper is as follows. Firstly in
Section~\ref{sec:main},
we set up a new CLT for LSSs of $B_n$ under elliptical distributions satisfying \eqref{BZ-condition}.
Then in  Section~\ref{sec:sph}, based on the derived results, we theoretically investigate the problem of sphericity test for covariance matrices. This is done by discussing a John's-type test from \cite{Tian15} for general alternative models and a likelihood ratio test from \cite{Onatski13} for {\em spiked covariances} under arbitrary elliptical distributions.
For illustration, the John-type test is applied to analyze a dataset of weekly stock returns in Section \ref{sec4}.  Technical proofs of the main results are gathered in Section~\ref{sec:proofs}. Some supporting lemmas are postponed to Appendix. The paper has also an on-line supplementary file which includes the following materials: (i) CLT for general moment LSSs; (ii) simulations regarding the John-type test; (iii) Proofs of some lemmas.

\section{High-dimensional theory for eigenvalues of $B_n$}
\label{sec:main}
This section investigates asymptotic behaviors of the eigenvalues of $B_n$, referred as {\em sample eigenvalues}. We begin with proposing an equivalent condition of \eqref{BZ-condition} under the settings of the elliptical model in \eqref{eds}.


\begin{lemma}\label{lem1}
	Suppose that a $p$-dimensional random vector $\x$ has a stochastic form $\x=\xi A\bu$ as defined in \eqref{eds} with the radius $\xi$ normalized as  $E(\xi^2)=p$. If the spectral norm of $\Sigma=AA'$ is uniformly bounded in $p$, then the following two conditions are equivalent:
	\begin{equation*}\label{mom4}
	{\rm a)}\quad  Var(\x' C_p\x)=o(p^2),\quad{\rm b)}\quad E(\xi^4)=p^2+o(p^2),
	\end{equation*}
	as $p\to\infty$, where $\{C_p\}$ is any sequence of symmetric matrices
	with bounded spectral norm.
\end{lemma}
\begin{remark}
	The fourth moment condition b) together with the normalization $E(\xi^2)=p$ characterize the class of elliptical distributions discussed in this paper.
	For the normal case, the squared radius $\xi^2\sim \chi_p^2$ and thus $E(\xi^2)=p$ and $E(\xi^4)=p^2+2p$.
	In general, the typical order of $E(\xi^4)$ is $p^2+\tau p+o(p)$ with $\tau\geq 0$ a constant. Hence a specific elliptical distribution can be recognized by evaluating the ratio
	\begin{equation}\label{mm-con}
	{E(\xi^4)}/{E^2(\xi^2)}=1+{\tau}/{p}+o(p^{-1}),
	\end{equation}
	as $p\rightarrow\infty.$
	We note that the parameter $\tau$ has a non-negligible effect on the limiting distributions of LSSs of $B_n$, see Theorem \ref{clt}.
		The  proof of Lemma \ref{lem1} is given in the supplementary material \citep{suppm}.
\end{remark}

In the following we provide three examples of elliptical family satisfying the condition \eqref{mm-con}.
Some commonly seen elliptical distributions are also checked and the results are summarized in Table \ref{table:s0}.
\begin{example}
A $p$-dimensional centered multivariate Pearson type II distribution has a density function
\begin{equation}
f(\x)=c_p|\Sigma_p|^{-\frac{1}{2}}\left[1-\x'\Sigma_p^{-1}\x\right]^{\frac{\beta}{2}-1},
\end{equation}
where $c_p=\pi^{-p/2}\Gamma[(\beta+p)/2]/\Gamma(\beta/2)$ and $\beta>0$. The stochastic representation of such a distribution is $\x=\xi \Sigma_p^{1/2} \bu$, where $\xi^2$ follows the distribution $Beta(p/2,\beta/2)$, see \cite{FZ90}. Therefore, we have
$$
{E(\xi^4)}/{E^2(\xi^2)}=1+{2\beta}/{(p^2+\beta p+2p)},
$$
which verifies the condition in \eqref{mm-con} with $\tau=0$.
\end{example}

\begin{example}
The family of Kotz-type distributions introduced by \cite{K75}  is an important class of elliptical distributions, which includes normal distributions, exponential power distributions, and double exponential distribution as special cases. The density function of a  centered  Kotz-type random variable $\x$ is
\begin{equation}
f(\x)=c_p|\Sigma_p|^{-\frac{1}{2}}\left[\x'\Sigma_p^{-1}\x\right]^{k-1}\exp\left\{-\beta\left[\x'\Sigma_p^{-1}\x\right]^s\right\},
\end{equation}
where $c_p=s\beta^{\alpha}\pi^{-p/2}\Gamma(p/2)\Gamma(\alpha)$ with $\alpha=(k-1+p/2)/s>0$ and $(\beta, s)>0$. Write $\x=\xi \Sigma_p^{1/2}\bu$. The $2s$ power of the radius is $\xi^{2s}=[\x'\Sigma_p^{-1}\x]^s$ which has the characteristic function
\begin{equation}\label{cf}
E(e^{it\xi^{2s}})=c_p\int e^{it[\x'\Sigma_p^{-1}\x]^{s}}f(\x)d\x\propto \int e^{itx}x^{\alpha-1}e^{-\beta x}dx,
\end{equation}
where the seconded integral is derived by polar coordinates transformation. This characteristic function implies that $\xi^{2s}$ follows the Gamma distribution $Gamma(\alpha, \beta)$. Simple calculations reveal that
$$
\frac{E(\xi^4)}{E^2(\xi^2)}=\frac{\Gamma(\alpha+2/s)\Gamma(\alpha)}{\Gamma^2(\alpha+1/s)}=1+\frac{1}{s^2\alpha}+o(\alpha^{-1}),
$$
which verifies the condition in \eqref{mm-con} with $\tau=2/s$.
For the mentioned three special cases, their details are presented in the 2-4th rows of Table \ref{table:s0}.
\end{example}
	
\begin{example}
Let $\x=\xi A\bu$ with $\xi^2=\sum_{j=1}^py_j^2$ independent of $\bu$, where $(y_j)$ is a sequence of i.i.d.\ random variables with
$$
E(y_1^2)=1\quad\text{and}\quad E(y_1^4)=\mu_4<\infty.
$$
Then it is simple to check that $E(\xi^2)=p$ and $E(\xi^4)=p^2+(\mu_4-1) p$ which verifies the condition in \eqref{mm-con} with $\tau=\mu_4-1$.
\end{example}


We should note that the condition \eqref{mom4} excludes some elliptical distributions, such as multivariate student-$t$ distributions and normal scale mixtures, as shown in the 5-6th rows of Table \ref{table:s0}. Indeed, sample eigenvalues from these distributions do not obey the generalized Mar\v{c}enko-Pastur law \citep{Karoui09, LY16}, which are then out of the scope of this paper.

\begin{table}[h]
	\setlength\tabcolsep{4pt}
	\begin{center}
		\caption{Some elliptical distributions and their verification of the condition \eqref{mom4}. The notation ``${\perp\!\!\!\perp}$" in the last row denotes independence.}
		\begin{tabular}{llllll}
			\hline
$\x=\xi A\bu\in\mathbb R^p$&Distribution of $\xi$&$E(\xi^4)/E^2(\xi^2)$&Condition \eqref{mom4} \\
\hline
Normal& $\xi^2\sim\chi_p^2$&$1+\frac{2}{p}$&Holds ($\tau=2$).\\
Double exponential& $\xi\sim Gamma(p,1)$&$1+\frac{4p+6}{p(p+1)}$&Holds ($\tau=4$).\\
Exponential power& $\xi^{2s}\sim Gamma\left(\frac{p}{2s},\frac{1}{2}\right)$&$1+\frac{2}{sp}+o(p^{-1})$&Holds ($\tau=\frac{2}{s}$).\\
Student-$t$&$\xi^2/p\sim F(p,v), v>4$&$1+\frac{2}{v-4}+\frac{2(v-2)}{p(v-4)}$& Not hold.\\
Normal scale mixture&$\xi^2=w\cdot v, w{\perp\!\!\!\perp}v, v\sim\chi_p^2$&$1+\frac{Var(w)}{E^2(w)}+\frac{2}{p}\frac{E(w^2)}{E^2(w)}$& Not hold.\\
			\hline
		\end{tabular}
		\label{table:s0}
	\end{center}
\end{table}

Now we are ready to investigate the asymptotic properties of sample eigenvalues in high-dimensional frameworks,
under the following assumptions.

\medskip
\noindent{\em Assumption}   (a). \quad
Both the sample size $n$ and dimension $p$ tend to infinity
in such a way that $c_n:=p/n\to c \in(0,\infty)$.

\medskip
\noindent{\em Assumption}   (b). \quad
There are two independent arrays of i.i.d.\ random variables
$(\bu_j)_{j\geq 1}$, $\bu_1\sim U(S^{p-1})$, and $(\xi_{j})_{j\geq 1}$ satisfying for some $\tau\geq 0$ and $\varepsilon>0$,
\begin{equation}\label{m-con}
E(\xi_1^2)=p,\quad E(\xi_1^4)=p^2+\tau p+o(p)\  \text{and}\ E\Big|\frac{\xi_1^2-p}{\sqrt{p}}\Big|^{2+\varepsilon}<\infty,
\end{equation}
such that  for each $p$ and $n$  the observation vectors can be represented as
$\x_j = \xi_j A\bu_j,$
where $A$ is a $p\times p$ matrix. 

\medskip
\noindent{\em Assumption}   (c). \quad
The spectral distribution $H_p$ of the matrix $\Sigma:=AA'$ weakly converges to a probability distribution $H$, as $p\rightarrow\infty$, referred as {\em Population Spectral Distribution} (PSD).
Moreover, the spectral norm of the sequence $(\Sigma)$ is uniformly bounded in $p$.

In the sequel, for any function $G$ of bounded variation on the real line, its Stieltjes transform is defined by
\begin{align}\label{Stie}
m(z)=\int\frac{1}{\lambda-z}dG(\lambda),~~z\in\mathbb{C}\setminus S_G,
\end{align}
where $S_G$ stands for the support of $G$.
Then we have the following theorems.
\begin{theorem}\label{lsd}
	Suppose that Assumptions (a)-(c) hold. Then, almost surely, the empirical spectral distribution $F^{B_n}$ converges weakly to a probability distribution $F^{c,H}$,
	whose Stieltjes transform $m=m(z)$ is the only  solution to the equation
	\begin{eqnarray}\label{mp1}
	m=\int\frac{1}{t(1-c-czm)-z}dH(t)~,\quad z\in\CC^+,
	\end{eqnarray}
	in the set $\{m\in \mathbb C:  -(1-c)/z+cm\in{\mathbb C^+}\}$ where $\mathbb C^+\equiv\{z\in \mathbb C: \Im(z)>0\}$.
\end{theorem}
\begin{remark}
	Theorem \ref{lsd} follows from Lemma \ref{lem1} and Theorem 1.1 in \cite{BZ08}, and thus we omit its proof here.
	Let $\um=\um(z)$ be the Stieltjes transform of $\underline{F}^{c,H}=cF^{c,H} + (1-c)\de_0$. Then Equation \eqref{mp1} can be re-expressed as
	\begin{equation}  \label{mp}
	z  =  - \frac1 {\um}  +  c \int\!\frac{t}{1+t\um} dH(t)~,\quad z\in\CC^+,
	\end{equation}
    which is the so-called Silverstein equation \citep{S95}.
\end{remark}

Let $F^{c_n, H_p}$ be the distribution defined by \eqref{mp1} with the parameters $(c,H)$ replaced by $(c_n, H_p)$ and denote $G_n=F^{B_n}-F^{c_n, H_p}$. We next study the fluctuation of centralized LSSs with form
\begin{eqnarray*}
	\int f(x) dG_n(x)=\int f(x) d[F^{B_n}(x)-F^{c_n,H_p}(x)],
\end{eqnarray*}
where $f$ is a function on the real line.

\begin{theorem}\label{clt}
	Suppose that Assumptions (a)-(c) hold.  Let $f_1,\ldots, f_k$ be functions analytic on an open interval containing
	\begin{equation}\label{interval}
	\left[\liminf_{p\rightarrow\infty}\lambda_{\min}^{\Sigma}\delta_{(0,1)}(c)(1-\sqrt{c})^2,\limsup_{p\rightarrow\infty}\lambda_{\max}^{\Sigma}(1+\sqrt{c})^2\right].
	\end{equation}
	Then the random vector
	$$
	p\left(\int f_1(x)dG_n(x),\ldots, \int f_k(x)dG_n(x)\right)
	$$
	converges weakly to a Gaussian vector $(X_{f_1},\ldots,X_{f_k})$, with mean function
	\begin{align*}
		 EX_f=&-\frac{1}{2\pi\rm i}\oint_{\mathcal C_1} f(z) \int \frac{c(\um'(z)t)^2}{\um(z)(1+t\um(z))^3}dH(t)dz\\
		&-\frac{\tau-2}{2\pi\rm i}\oint_{\mathcal C_1} f(z) \int\frac{(z\um(z)+1)\um'(z)t}{(1+t\um(z))^2}dH(t)dz
	\end{align*}
	and covariance function
	\begin{align*}
		{ Cov}\left(X_f,X_g\right)
		=&-\frac{1}{2\pi^2}\oint_{\mathcal C_1}\oint_{\mathcal C_2}\frac{f(z_1)g(z_2)\um'(z_1)\um'(z_2)}{(\um(z_1)-\um(z_2))^2}dz_1dz_2\nonumber\\
		&+c(\tau-2)\int xf'(x)dF(x)\int xg'(x)dF(x),\label{lss-cov}
	\end{align*}
	$(f, g \in \{f_1,\cdots,f_k\})$, where the contours $\mathcal C_1$ and $\mathcal C_2$ are non-overlapping, closed, counter-clockwise orientated in the complex plane, and each enclosing the support of
	the limiting spectral distribution $F^{c,H}$.
\end{theorem}

\begin{remark}\label{rem2.3}
 When the population is normal, or rather  $\tau=2$,  this theorem coincides with the main result in \cite{BS04}. It implies that the high order correlation among the components of the population affects both the limiting mean vectors and the covariance matrices of LSSs by additive quantities proportional to $\tau-2$. This factor can be further decomposed into two parts, $\tau$ and $-2$, which correspond respectively to the effect from the radius $\xi$ and that from the direction $\bu$ (considering the case $\xi^2\equiv p$). It's interesting to see that these two kinds of dependency have opposite effects and they may cancel each other for normal population.
\end{remark}

As an application of Theorem \ref{clt}, we consider $\hat \beta_{nj}=\int x^jdF^{B_n}(x), \quad j=1,2$, the first two moments of sample eigenvalues. Theorem \ref{clt} implies
\begin{eqnarray}\label{clt-beta}
p\left(
\begin{matrix}
\hat \beta_{n1}-  \beta_{n1}\\
\hat \beta_{n2}-  \beta_{n2}
\end{matrix}\right)\xrightarrow{D}
N\left(\left(
\begin{matrix}
v_1\\
v_2
\end{matrix}
\right),
\left(\begin{matrix}
\psi_{11}&\psi_{12}\\
\psi_{12}&\psi_{22}
\end{matrix}
\right)\right),
\end{eqnarray}
where the parameters possess explicit expressions as
\begin{align*}
	&\beta_{n1}= \gamma_{n1},\quad  \beta_{n2}=\gamma_{n2}+ c_n\gamma_{n1}^2,\quad v_1=0,\quad v_2=c\gamma_2+c(\tau-2)\gamma_1,\\
	&\psi_{11}=2 c \gamma_2+c (\tau-2) \gamma_1^2,\quad
	\psi_{12}=4c\gamma_3+4 c^2 \gamma_1 \gamma_2 + 2c(\tau-2) \gamma_1 (c \gamma_1^2 + \gamma_2),\\
	&\psi_{22}=8 c\gamma_4+4c^2\gamma_2^2  + 16c^2\gamma_1 \gamma_3 + 8c^3\gamma_1^2 \gamma_2 +4c(\tau-2) (c \gamma_1^2 + \gamma_2)^2,
\end{align*}
where $\gamma_{nj}=\int t^jd H_p(t)$ and $\gamma_{j}=\int t^jd H(t)$ for $j\geq1$.
For LSSs of higher order moments, explicit formulas of their limiting means and covairances are discussed in the supplementary material \citep{suppm}.
\begin{figure}[htbp]
	\begin{minipage}[t]{0.5\linewidth}
		\includegraphics[width=2.6in]{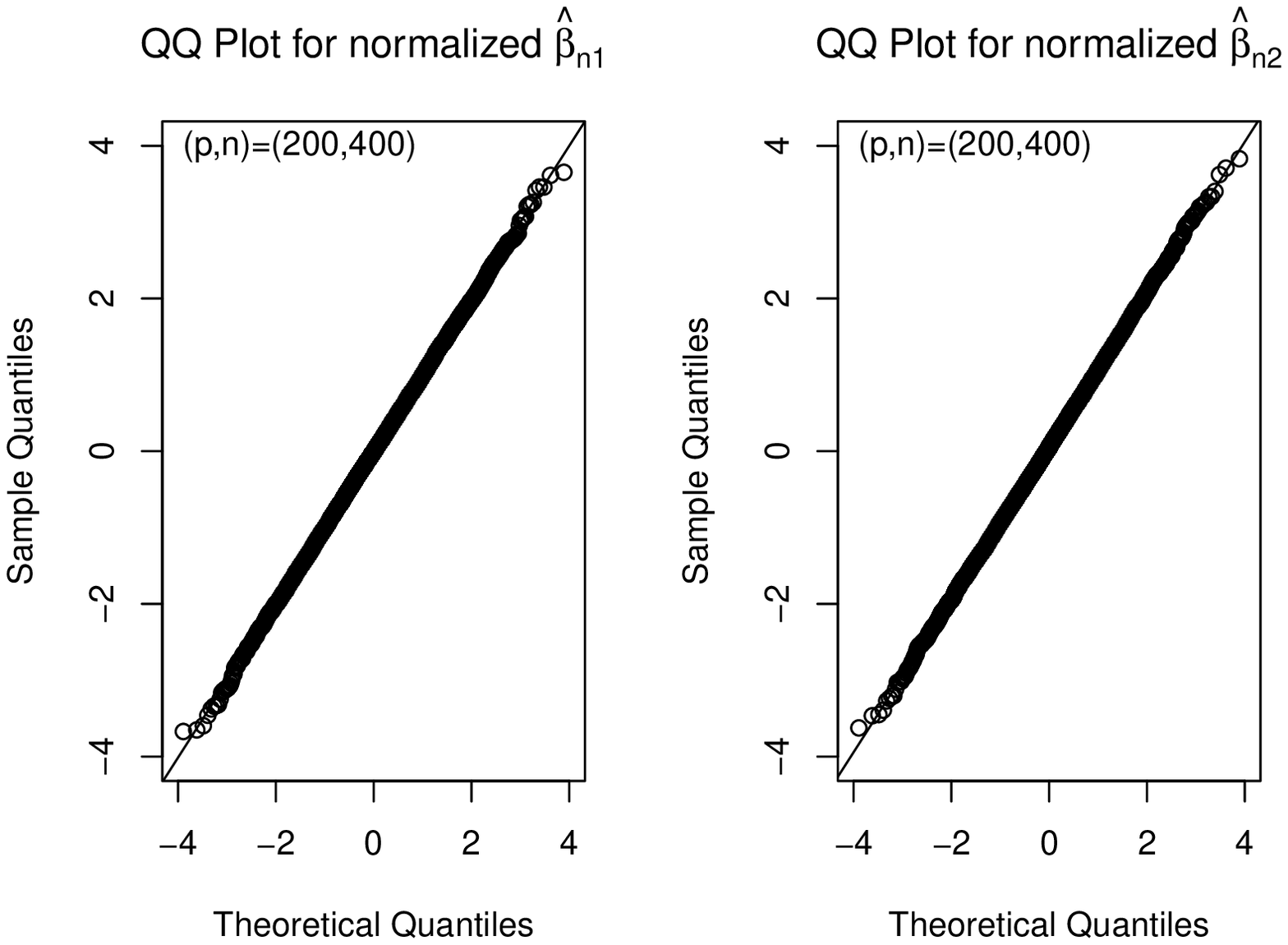}
	\end{minipage}%
	\begin{minipage}[t]{0.5\linewidth}
		\includegraphics[width=2.6in]{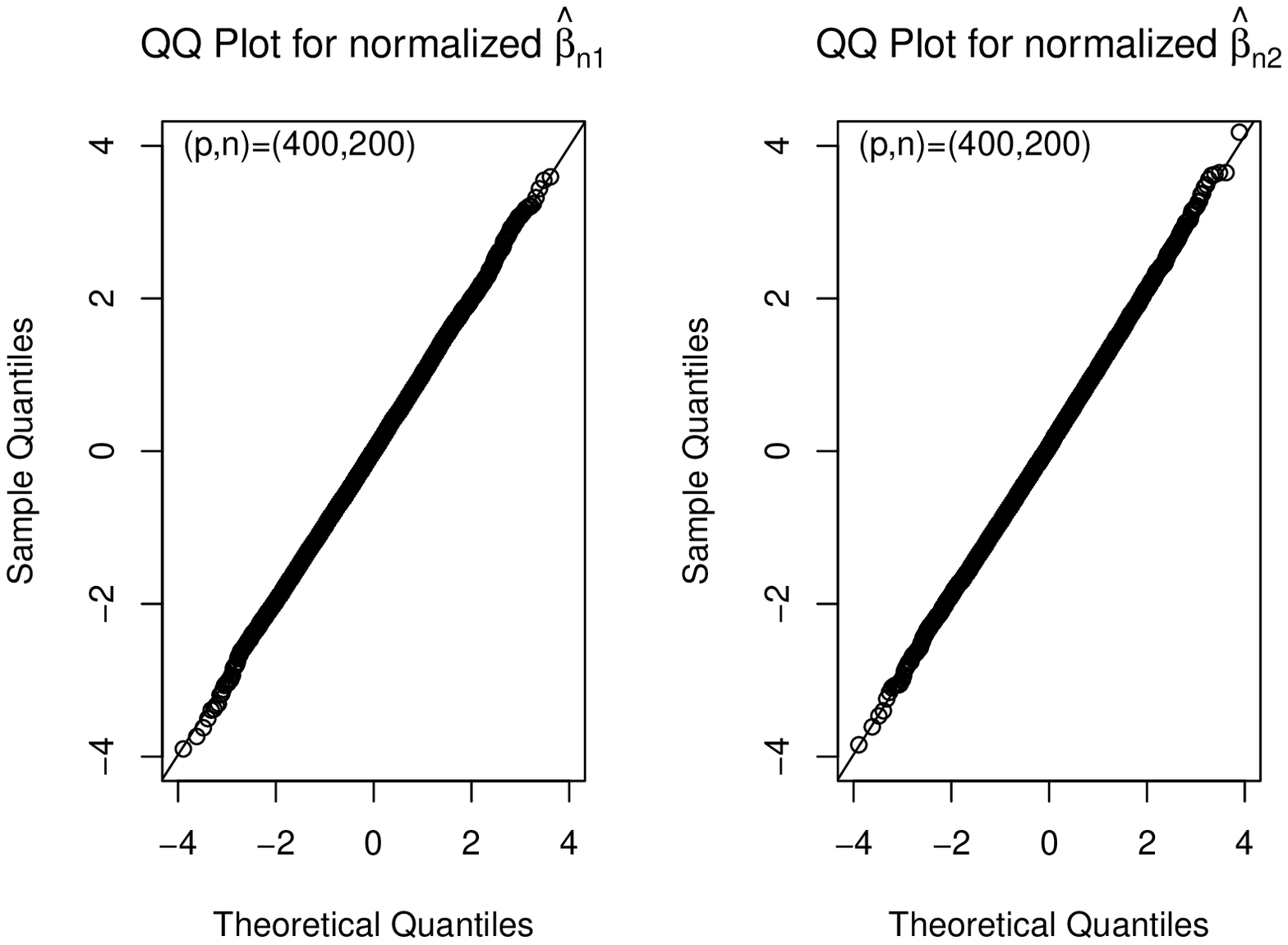}
	\end{minipage}
	\begin{minipage}[t]{0.5\linewidth}
		\includegraphics[width=2.6in]{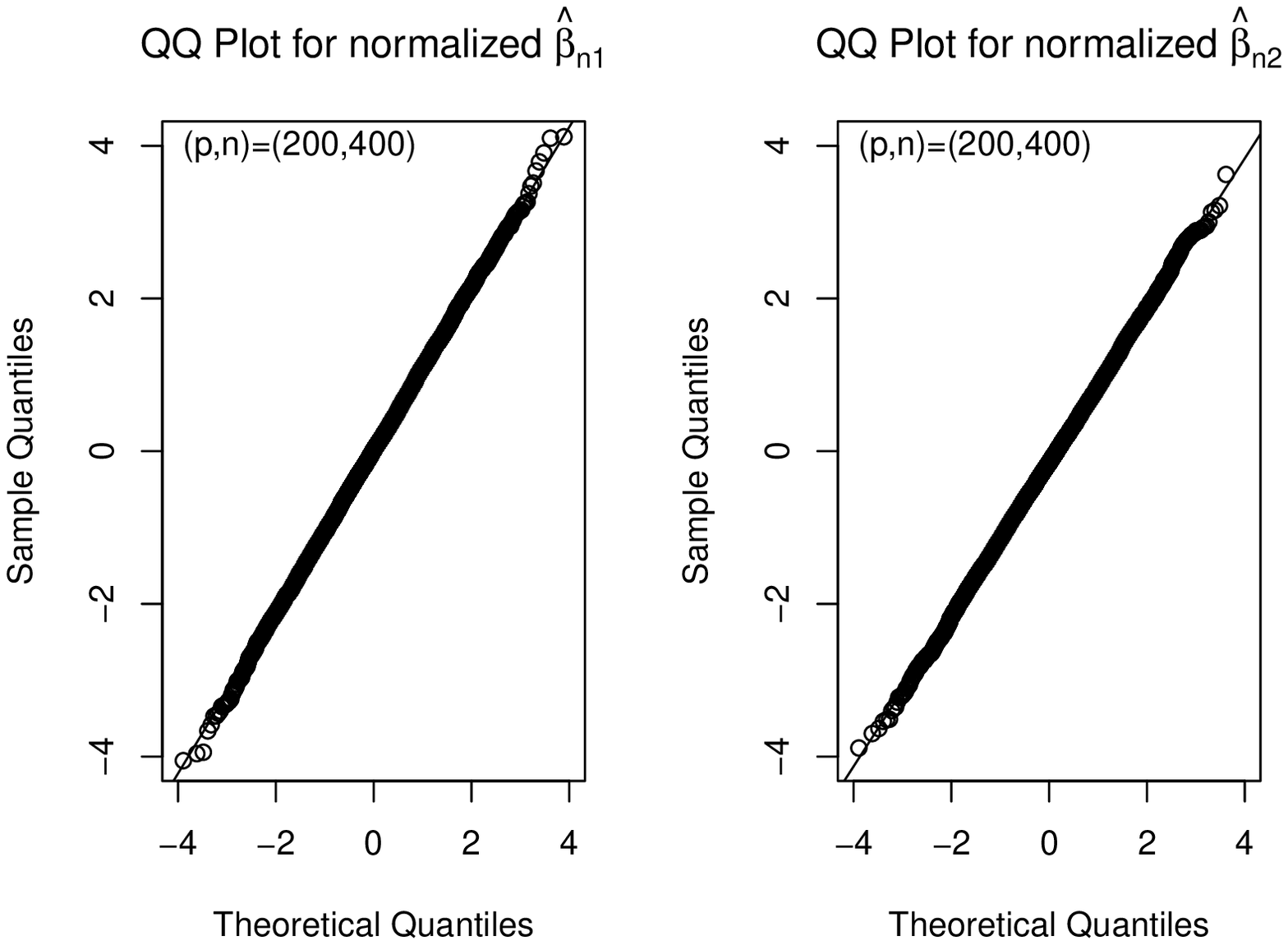}
	\end{minipage}%
	\begin{minipage}[t]{0.5\linewidth}
		\includegraphics[width=2.6in]{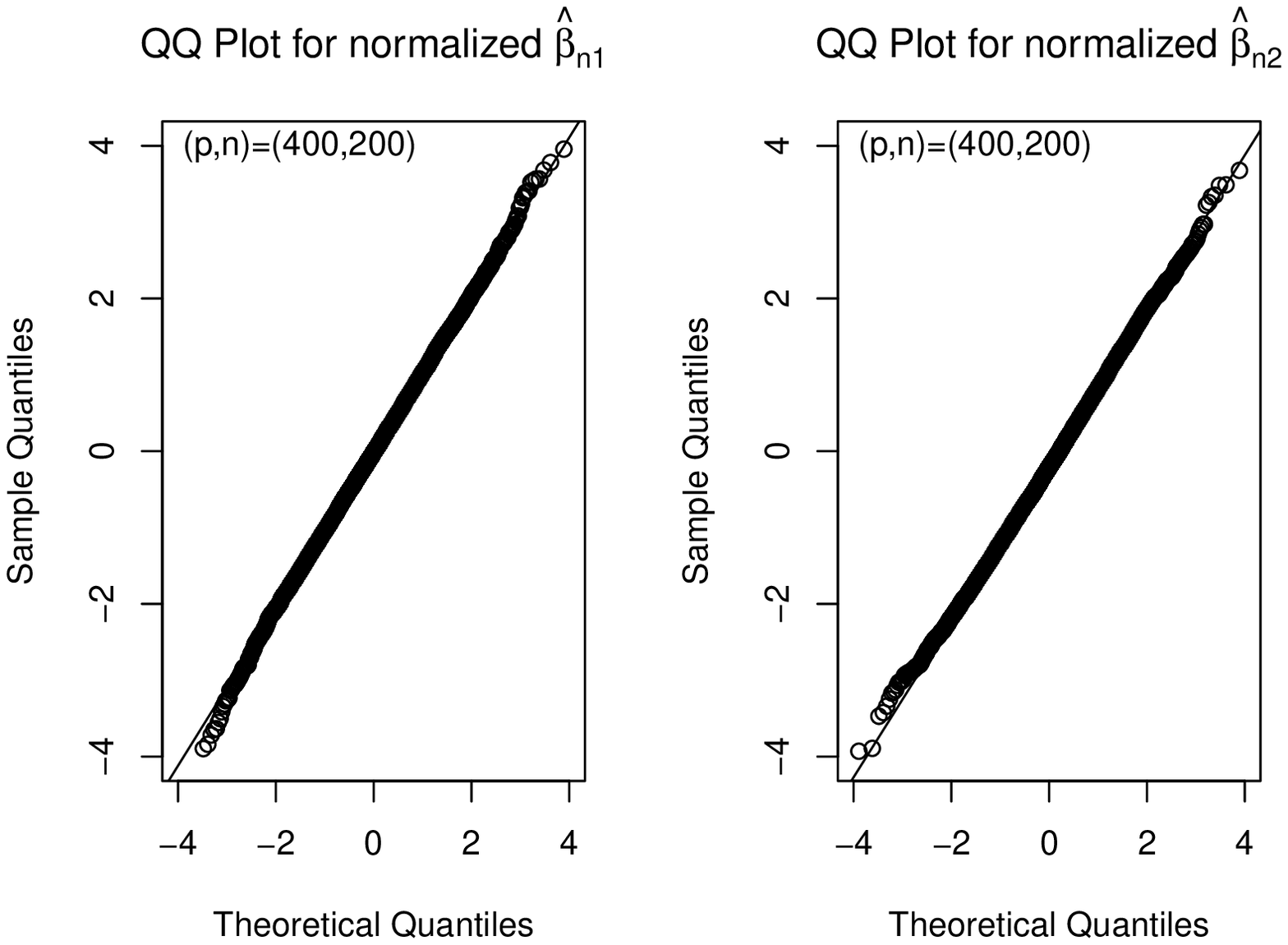}
	\end{minipage}
	\caption{Normal QQ-plots for normalized $\hat\beta_{n1}$ and $\hat{\beta}_{n2}$ from 10,000 independent replications.
		Upper panels: $\xi\sim k_1Gamma(p,1)$ with $k_1=1/\sqrt{p+1}$. Lower panels: $\xi^2\sim k_2Beta(p/2,2)$ with $k_2=p+4$. The dimensional settings are $(p,n,c)=(200,400,0.5), (400,200,2)$. }
	\label{fig0}
\end{figure}

We conduct a small simulation experiment to examine the fluctuations of $\hat\beta_{n1}$ and  $\hat\beta_{n2}$. In the experiment, the PSD $H_p$ is fixed at $H_p=0.5\delta_1+0.5\delta_2$. The distribution of $\xi$ is  selected as (1) $\xi\sim k_1Gamma(p,1)$ with $k_1=1/\sqrt{p+1}$ and (2) $\xi^2\sim k_2Beta(p/2,2)$ with $k_2=p+4$, which correspond the CLT with $\tau=4$ and $\tau=0$, respectively. The factors $k_1$ and $k_2$ are selected to satisfy $E(\xi^2)=p$. The dimensional settings are $(p,n,c)=(200,400,0.5), (400,200,2)$ and the number of independent replications is $10,000$. Normal QQ-plots for normalized statistics, i.e. $p(\hat{\beta}_{n1}-\beta_{n1})/\sqrt{\psi_{11}}$ and $[p(\hat{\beta}_{n1}-\beta_{n1})-v_2]/\sqrt{\psi_{22}}$, are displayed in Figure 1. Their asymptotic standard normality
is well confirmed in all studied cases.

\section{Testing for high-dimensional spherical distributions}
\label{sec:sph}

\subsection{John's test and its extension}

In this section,  we revisit the sphericity test for covariance matrices in high-dimensional frameworks. For this particular test probelm, the underlying population can follow arbitrary elliptical distribution, which may violate the condition in \eqref{BZ-condition}.

The sphericity test on the covariance matrix $\Sigma$ is
\begin{equation}\label{hp}
H_{0}:\Sigma= \sigma^2 I_p \quad\text{v.s.}\quad H_{1}:\Sigma \neq \sigma^2 I_p,
\end{equation}
where $\sigma^2$ is an unknown scalar parameter.
When the dimension $p$ is fixed, for normal populations,  \cite{John(1972)} proposed a locally most powerful invariant test statistic to deal with the sphericity hypothesis based on the spectrum of sample covariance matrices.
Due to its concise form and broad applicability, this kind of  test is quite favorable for high dimensional situations and has been extensively studied in recent years. See, for example, \cite{LW02}, \cite{WangYao13}, \cite{Tian15} for the linear transform model in \eqref{linear}, while \cite{Z14} and \cite{PV16} for the elliptical model in \eqref{eds}.
In particular, the test statistic in \cite{Tian15}  synthesizes the first four moments of sample eigenvalues, by which it gains extra powers for spike-like alternative covariance matrices.
However this statistic is not valid for general elliptical populations \citep{LY16}.
Hence, we next develop an analogue test procedure with the help of the theoretical results in Section \ref{sec:main}, and then compare it numerically with that from \cite{PV16}.

Since the hypotheses in \eqref{hp} are only concerned with the shape component of $\Sigma$,  by convention, we transform the original samples into the so-called spatial-sign samples,  that is,
$$\check\x_j := \sqrt{p}\x_j/||\x_j||=\sqrt{p}A\bu_j/||A\bu_j||,\ j=1,\ldots,n.$$
Therefore, testing the sphericity of $\Sigma$ can now be converted to testing the identity of $\check\Sigma=E(\check \x_1\check\x_1')$. This inference can be realized by constructing spectral statistics of
$\check B_n=\sum_{j=1}^n\check\x_j\check\x_j'/n$.
Specifically, let
$$
\alpha_{nj}=p^{-1}\tr(\check\Sigma^j)\quad\text{and}\quad \check\beta_{nj}=p^{-1}\tr(\check B_n^j),
$$
$j=0,1,2,\ldots.$
By verifying the condition in \eqref{BZ-condition} for $\check\x_1$, 
one may conclude that Theorem \ref{lsd} also holds for $(\check\Sigma, \check B_n)$ with all conditions on $\xi$ removed.
Then, similar to \cite{Tian15}, from the fact that $\check\beta_{n1}\equiv1$, one may obtain estimators of $\alpha_{n2}$ and $\alpha_{n4}$ as
\begin{eqnarray*}
	\check\alpha_{n2}=\check\beta_{n2}-c_n,\quad \check\alpha_{n4}=\check\beta_{n4}-4c_n\check\beta_{n3}-2c_n(\check\beta_{n2})^2+10c_n^2\check\beta_{n2}-5c_n^3,
\end{eqnarray*}
respectively, and two simple statistics for the sphericity test as
$$
T_1=\check\alpha_{n2}-1\quad\text{and}\quad T_2=\check\alpha_{n4}-1.
$$
Moreover, their joint null distribution is directly from \eqref{clt-beta} with $\tau=0$.
\begin{theorem}\label{th3}
	Suppose that Assumptions {\rm (a)--(c)} [removing the moment conditions in \eqref{m-con}] hold. Under the null hypothesis,
	\begin{equation*}
	n\left(T_1,T_2\right)
	\xrightarrow{D}N_2(\mu,\Omega),
	\end{equation*}
	where $\mu=(-1, -6+c)$ and the covariance matrix $\Omega=(\omega_{ij})$ with $\omega_{11}=4,\omega_{12}=24,$ and $\omega_{22}=8(18 + 12 c +  c^2).$
\end{theorem}

The two statistics $T_1$ and $T_2$, together with their null distributions, provide two test procedures for the identity of $\check\Sigma$ (thus the sphericity of $\Sigma$).
The test statistic $T_1$ agrees with that in \cite{PV16}, where its null asymptotic distribution is proved to be universal whenever $\min\{n,p\}\to\infty$. For the case where the population mean is unknown, see \cite{Z14}.
The test statistic $T_2$ is new. Compared with $T_1$, it is more sensitive to  extreme eigenvalues of $\check\Sigma$ and thus can serve as a complement of $T_1$.
Parallel to \cite{Tian15}, a joint statistic of $T_1$ and $T_2$ can be constructed as
\begin{equation*}\label{tm}
T_{m}=\max\left\{\frac{nT_1+1}{2},\frac{nT_2+6-c_n}{\sqrt{8(18+12c_n+c_n^2)}}\right\},
\end{equation*}
where the two original statistics are both standardized according to their asymptotic null distributions.

\begin{theorem}\label{th4}
	Suppose that Assumptions {\rm (a)--(c)} [removing the moment conditions in \eqref{m-con}] hold and let $\delta_p=p\tr(\Sigma^2)/\tr^2(\Sigma)-1$.
	\begin{itemize}
		\item[\rm(i)] Under the null hypothesis, for any $x\in \mathbb R$,
		\begin{eqnarray*}
			P\left(T_m\leq x\right) \rightarrow
			\int_{-\infty}^x\int_{-\infty}^x\frac{1}{2\pi\sqrt{1-\rho^2}}\exp\left\{-\frac{u^2-2\rho uv+v^2}{2(1-\rho^2)}\right\}dudv,
		\end{eqnarray*}
		where $\rho=6/\sqrt{2(18+12c+c^2)}$.
		
		\item[\rm (ii)]
		Under the alternative hypothesis, if $n\delta_p\to\infty$ then the power of the test $T_m$ goes to 1 as $(n,p)\to \infty$.
	\end{itemize}
\end{theorem}

The asymptotic null distribution of $T_m$ is an immediate consequence of Theorem \ref{th3}. The consistency of $T_m$ can be proved by showing either the consistency of $T_1$ or  $T_2$.
As the consistency of $T_1$ has been given in \cite{Z14}, we omit its proof.

We have run a simulation experiment for the tests $T_1$, $T_2$, and $T_m$ to check their finite-sample properties under similar model settings as in \cite{Tian15}. The results show that all the three tests have satisfactory empirical sizes and powers. In addition, compared with $T_1$ and $T_2$, the test $T_m$ exhibits its robustness against different types of alternative models, see the supplementary material \cite{suppm}.

\subsection{Sphericity test under spiked alternative model}

The sphericity test $T_m$ applies to general alternative models. However, its consistency requires $n\delta_p\to\infty$ which excludes the well-known {\em spiked covariance} model \citep{Johnston01}. For the simplest spiked model, the covariance matrix can be expressed as $\Sigma=\sigma^2(I_p+h vv')$ where $\sigma^2$ and $I_p$ are as before, $h$ is a constant, and $v$ is a unit vector in $\mathbb R^p$. Both $h$ and $v$ are unknown parameters. Thus the sphericity hypotheses in \eqref{hp} reduce to
\begin{equation}\label{hps}
H_0: h=0\quad v.s.\quad H_1: h>0.
\end{equation}
It's obvious that $T_m$ will asymptotically fail to reject  such alternatives since $n\delta_p\to 0$.  What's more, this testing problem will become more difficult but attractive when the signal $h$ falls below the threshold $\sqrt{c}$,
see \cite{BR13,Onatski13,Onatski14,D15}, and references therein. Hence, applying the CLT for LSSs under elliptical distributions, we discuss a test procedure for \eqref{hps} proposed by \cite{Onatski13}, which was built under normal populations.

In \cite{Onatski13}, the authors discussed a likelihood ratio test based on the joint distribution of sample eigenvalues from normal populations. This test was especially designed for the local alternative $H_1: h\in (0,\sqrt{c})$ and the employed statistic was approximated by a special LSS. In our settings, this LSS can be formulated as
\begin{equation}
T_{LR}(s)=\int \ln(z(s)-x)dF^{\check B_n}(z)-\int \ln(z(s)-x)dF^{c_n, \delta_1}(x),
\end{equation}
where $s\in (0, \bar s)$ is a testing parameter and $z(s)=(1+s)(c_n+s)/s$.  The upper bound $\bar s$ of $s$ is chosen as $\bar s=\sqrt{c_n}$ for $h\in[0,\sqrt{c_n}~]$ and $\bar s=h^{-1}c_n$ for $h\in(\sqrt{c_n}, \infty)$ such that $z(s)$ is larger than the limit of $\lambda_{\max}(\check B_n)$, the largest sample eigenvalues.
Applying Theorem \ref{clt}, one may get the asymptotic distribution of $T_{LR}(s)$ under general elliptical distributions.

\begin{theorem}\label{clt-spike}
	Suppose that Assumptions {\rm (a)--(c)} [removing the moment conditions in \eqref{m-con}] hold. Under the null hypothesis, for any fixed $s\in(0, \bar{s})$,
	\begin{align}
	pT_{LR}(s)\xrightarrow{D}N(\mu_{s},\sigma_s^2),
	\end{align}
	where the respective mean and variance functions are 	
	\begin{align*}
	\mu_{s}=\frac{1}{2}\ln(1-c^{-1}s^2)+c^{-1}s^2\ \text{and}\
	\sigma_s^2=-2\ln(1-c^{-1}s^2)-2c^{-1}s^2.
	\end{align*}
\end{theorem}
The  proof of Theorem \ref{clt-spike} is given in the supplementary material \citep{suppm}.
Given a value of $s\in(0, \bar s)$  and a significance level $\alpha\in(0,1)$, the test $T_{LW}(s)$ rejects $H_0$ if $pT_{LW}(s)<\sigma_s\Phi^{-1}(\alpha)+\mu_{s}$,
where $\Phi(x)$ denotes the standard normal distribution function.  Unlike \cite{Onatski13}, the theoretical power of this test is not available at present since $pT_{LR}(s)$ is not a likelihood ratio statistic in elliptical distributions. Another reason is that Theorem \ref{clt} is inapplicable to this situation since the spatial-sign sample is not anymore elliptically distributed under $H_1$.

Let's take a step back and only consider the testing problem in elliptical distributions satisfying \eqref{mm-con}.
For simplicity, we assume $\sigma$ is known and set $\sigma =1$, so that the test $T_{LR}(s)$ is still valid by simply substituting the sample covairance $B_n$ into $\check B_n$, i.e.,
\begin{align*}
\widetilde T_{LR}(s)=\int \ln(z(s)-x)dF^{ B_n}(z)-\int \ln(z(s)-x)dF^{c_n, \delta_1}(x),
\end{align*}
whose asymptotic distribution under both the null and alternative hypotheses is described in the following theorem.

\begin{theorem}\label{clt-spike-2}
	Suppose that Assumptions {\rm (a)--(c)} hold. Let $h_0$ be the true value of $h$ and $\sigma=1$, then for any fixed $s\in(0, \bar{s})$,
	\begin{align}
	p\widetilde T_{LR}(s)\xrightarrow{D}N(\tilde\mu_{s}(h_0),\tilde\sigma_s^2),
	\end{align}
	where the respective mean and variance functions are
	\begin{align*}
	\tilde\mu_{s}(h)&=\frac{1}{2}\ln(1-c^{-1}s^2)+(1-\tau/2)c^{-1}s^2+\ln(1-c^{-1}sh),\\
	\tilde\sigma_s^2&=-2\ln(1-c^{-1}s^2)-(2-\tau)c^{-1}s^2.
	\end{align*}
\end{theorem}
This theorem is a direct conclusion of  Theorem \ref{clt}. It proof is similar to that of Theorem \ref{clt-spike} and we thus omit it here. From Theorem \ref{clt-spike-2}, the power function of $\widetilde T_{LW}(s)$ is
$$
P_{H_1}(\widetilde T_{LW}(s)\ \text{reject}\ H_0)=\Phi\left[\Phi^{-1}(\alpha)-\frac{\tilde\mu_{s}(h_0)-\tilde \mu_{s}(0)}{\tilde\sigma_s}\right],\quad h_0>0.
$$
For normal populations ($\tau=2$), this power function reaches its maximum at $s=h_0$, which agrees with (5.1) in Proposition 9 of \cite{Onatski13}. In general, the maximizer may not locate at $h_0$. An interesting case is $\tau=0$, for which the power function tends to 1 as $s\to 0^{+}$. This is from the fact that
\begin{align*}
\tilde\mu_{s}(h_0)-\tilde\mu_{s}(0)=-c^{-1}sh_0+o(s)~\text{and}~
\sigma_s^2=2c^{-2}s^4+o(s^4).
\end{align*}
At this time,   $\widetilde T_{LW}(s)$ will successfully detect any  positive $h_0$ as long as $s$ is close to zero.


\section{An empirical study}\label{sec4}

For illustration, we apply the test procedure based on $T_m$ to analyze weekly returns of the stocks from S\&P 500. The tests $T_{LR}$ and $\tilde T_{LR}$ are not included in this analysis since there is a lack of evidence to fit the data using the simplest spiked model. According to The North American Industry Classification System (NAICS),  which  is used by business and government to classify business establishments, the 500 stocks can be divided into 20 sectors.  Nine of them are removed from our analysis since their numbers of stocks are all less than 10.  The remaining 11 sectors as well as their numbers of stocks are listed  in Table \ref{tb3}.
\begin{table}[h]
		\caption{Number of stocks in each NAICS Sectors. 
	}\label{tb3}
	\centering
	\begin{tabular}{c|c|c|c|c|c|c|c|c|c|c|c}
		\hline
		Sector&\  1\ \ &\  2\ \ &\ \ 3\ \  &\  4\  &\ \ 5\ \ & \ 6\ \ & \ 7\ \ & \ 8\ \ & \ 9\ \ & \ 10\ \ & \ 11\ \ \\\hline
		Number of stocks&30&32&189&17&36&14&37&65&14&17&11\\\hline
	\end{tabular}
\end{table}
Usually the stocks in the same sector are correlated, and the stocks in different sectors are uncorrelated. So it is expected that the weekly returns of stocks in the same sector are not spherically distributed, and it is interesting to see if the weekly returns of stocks in different sectors are spherically distributed. In the following, we apply $T_m$ to stocks in the same sector and stocks in different sectors respectively.

The original data are the  closing prices or the  bid/ask average of these stocks for the trading days in the first half of  2013, i.e., from 1 January 2013 to 30 June 2013, with total 124 trading days. This  dataset is downloadable from  the Center for Research in Security Prices  Daily Stock in Wharton Research Data Services.
The testing model is  established  as  follows.  Denote $p_l$ as the number of stocks  in the
$l$th sector, $u_{ij}(l)$ as the price of the $i$th stock  in the
$l$th sector on Wednesday of the $j$th week. The reason that we choose Wednesday's price here is   to avoid the weekend effect in  stock market. Thus we get  22 returns for each stock.    In order to  meet the condition of  the proposed test, the original data $u_{ij}(l)$ should be transformed by logarithmic difference, which is a very commonly used procedure in finance. There are a number of theoretical and practical advantages of using logarithmic returns. One of them is that the sequence of logarithmic returns are independent of each other for  big  time scales (e.g.  $\geq $ 1 day, see \cite{Rama(01)}).
Denote   $x_{ij}(l)=\ln(u_{i(j+1)}(l)/u_{ij}(l))$, $j=1,\dots, 21$ and $X(l)=(x_{ij}(l))_{p_l\times n}$, where $n=21$ is the   sample size.

Now  applying $T_m$ to the  dataset $X(l)$, $l=1,\dots,11$, respectively,  we obtain 11 $p$-values, which are all below $10^{-9}$. Therefore, we have strong evidence to believe that stocks in the same sector are not spherically distributed. This is consistent with our intuition. Next, we consider stocks in  different sectors. Specifically, we choose one stock from each sector to make up a group of 11 cross-sectoral stocks and then test  whether these stocks are spherically distributed.  Because there are about $9.79\times 10^{15}$ different groups, we just randomly draw 1,000,000 groups from them to analyze.  It turns out that the largest $p$-value is 0.3889, 231 $p$-values  are bigger than 0.05, and 69 $p$-values are bigger than 0.1.  These results again demonstrate that, when the number of stocks is not very small, it is hard to say weekly logarithmic returns for the
stocks are spherically distributed. It is also very interesting to analyze these spherically distributed stocks in different sectors, which have almost the same variances.

\section{Proof of Theorem \ref{clt}}\label{sec:proofs}

The proof of Theorem \ref{clt} relies on analyzing the resolvent of the sample covariance matrix $B_n$ and the general strategy follows the approach in \cite{BS04}. Also see \cite{BJPZ15} and \cite{Gao16} for recent developments.  However, as we are dealing with the new model equipped with nonlinear dependency, all technical steps of implementing this strategy have to be updated, or at least revalidated. They are presented in this section.
\subsection{ Sketch of the proof of Theorem \ref{clt}}

Let $v_0 > 0$ be arbitrary,  $x_r$ any number greater than the right end point of interval \eqref{interval}, and $x_l$ any negative number if the left end point of \eqref{interval} is zero, otherwise choose
$x_l\in(0,\liminf_{p\rightarrow\infty}\lambda_{\min}^{\Sigma}(1-\sqrt{c})^2)$. Let $ \mathcal C_u=\{x\pm iv_0: x\in[x_l,x_r]\}$ and define a contour $\mathcal C$
\begin{equation}\label{cont}
\mathcal C=\{x+iv: x\in \{x_r,x_l\}, v\in[-v_0,v_0]\}\cup \mathcal C_u.
\end{equation}
By definition, this contour encloses a rectangular region in the complex plane containing the support of the LSD $F^{c, H}$.
Denote by $m_n(z)$ and $m_{F^{c_n,H_p}}(z)$ the Stieltjes transforms of the ESD $F^{B_n}$ and the LSD $F^{c_n,H_p}$, respectively.
Their companion Stieltjes transforms are given by
$$\um_n(z)=-\frac{1-c_n}{z}+c_nm_n(z)\quad\text{and}\quad \um_{F^{c_n,H_p}}(z)=-\frac{1-c_n}{z}+c_nm_{F^{c_n,H_p}}(z).$$
With these notation, we define an empirical process on $\mathcal C$ as
$$
M_n(z)=p\left[m_n(z)-m_{F^{c_n,H_p}}(z)\right]=n\left[\um_n(z)-\um_{F^{c_n,H_p}}(z)\right],\quad z\in\mathcal C.
$$
Since $f_\ell,~\ell=1,\dots,k,$ in Theorem \ref{clt} are  analytic on an open region containing the interval \eqref{interval} (thus analytic on the region enclosed by $\mathcal C$), by Cauchy's integral formula, we have for any $k$ complex numbers $a_1,\dots,a_k$,
\begin{align*}
	\sum_{\ell= 1}^kpa_\ell\int f_\ell(x)dG_n(x)=-\sum_{\ell= 1}^k\frac{a_\ell}{2\pi i}\oint_{\mathcal C}f_\ell(z) M_n(z)dz,
\end{align*}
when all sample eigenvalues fall in the interval $(x_l, x_r)$, which is correct with overwhelming probability.
In order to remove the small probability event that some   sample eigenvalues fall outside the interval, we need a truncated version of $M_n(z)$, denoted by $\widehat{M}_n(z)$. Specifically,
let $\{\varepsilon_n\}$ be a sequence decreasing to zero satisfying  $\varepsilon_n>n^{-a}$ for some $a\in(3/4,1)$. The truncated process $\widehat{M}_n(z)$ for $z=x+iv\in\mathcal C$ is given by
\begin{eqnarray}
\widehat M_n(z)=
\begin{cases}
M_n(z)&z\in \mathcal C_n,\\
M_n(x+in^{-1}\varepsilon_n)& x\in \{x_l,x_r\},~ v\in[0, n^{-1}\varepsilon_n],\\
M_n(x-in^{-1}\varepsilon_n)& x\in \{x_l,x_r\},~ v\in[-n^{-1}\varepsilon_n, 0],
\end{cases}
\end{eqnarray}
where
$$\mathcal C_n=\mathcal C_u\cup \{x\pm iv: x\in\{x_l, x_r\}, v\in[n^{-1}\varepsilon_n, v_0]\},$$
on which $\widehat M_n(z)$ agrees with $M_n(z)$, is a regularized set of $\mathcal C$ excluding a small segment near the real line.
Then we have
\begin{lemma}\label{hatMn}
	Under the same assumptions in Theorem \ref{clt}, we have for any $\ell>0$,
	\begin{align*}
	\oint_{\mathcal C}f_\ell M_n(z)dz=\int_{\mathcal C_n} f_\ell(z)\widehat M_n(z)dz+o_p(1).
	\end{align*}
\end{lemma}
The  proof of  this lemma will  be  put in the supplementary
material, \cite{suppm}.
Hence, Theorem \ref{clt} follows by similar arguments on Pages 562-563 in \cite{BS04} and the following Lemma.

\begin{lemma}\label{mainlemma}
	Under Assumptions (a)-(c), the random process $\widehat M_n(\cdot)$ converges weakly to a two-dimensional Gaussian process $M(\cdot)$ with the mean function
	\begin{align}\label{mean}
	EM(z)=&\int\frac{c(\um'(z)t)^2dH(t)}{\um(z)(1+t\um(z))^3}+(\tau-2)\int\frac{(z\um(z)+1)\um'(z)tdH(t)}{(1+t\um(z))^2}
	\end{align}
	and covariance function
	\begin{align}
	Cov\left(M(z_1),M(z_2)\right)
	=&\frac{2\um'(z_1)\um'(z_2)}{(\um(z_1)-\um(z_2))^2}-\frac{2}{(z_1-z_2)^2}\nonumber\\
	&+\frac{\tau-2}{c}\left(\um(z_1)+z_1\um'(z_1)\right)\left(\um(z_2)+z_2\um'(z_2)\right),\label{cov-ij}
	\end{align}
where $z, z_1, z_2\in \mathcal C$.
\end{lemma}
The proof of this lemma is the main task of this section and  can be achieved by four steps as described below.
 Notice that in the proof we will use several inequalities frequently, which are presented as  lemmas in Appendix. We will show how and where to use these lemmas in the following.
Write for $z \in \mathcal C_n$,
\begin{align*}
\widehat M_n(z) &=p[m_n (z) - Em_n (z)]+p[Em_n (z)-m_{F^{ c_n,H_p}} (z)]\\
&:= M_{n}^{(1)}(z) + M_{n}^{(2)}(z).
\end{align*}
\begin{itemize}
	\item[] {\bf Step 1}: Truncation and rescaling of $\xi$.  This step regularizes the variables
    $\{\xi_j\}$ in $B_n= \sum_{j=1}^n\xi_j^2A\bu_j\bu_j'A'/n$ such that $\{\xi_j\}$ have proper bound for finite $(n,p)$ while maintaining the limiting distribution of the LSSs. Compared with the proof in \cite{BS04}, this result is entirely new since their model doesn't include  a radius variable at all. Moreover, our key inequalities (Lemmas \ref{ineq}-\ref{lambda-bound}) are all built on this regularization and thus their proofs have to be updated to accommodate the elliptical model.
	\item[] {\bf Step 2}: Finite dimensional convergence of $M_n^{(1)}(z)$ in distribution. This step finds the joint limiting distribution of an $r$-dimensional vector $\big(M^{(1)}_n(z_\ell)\big)_{1\leq\ell\leq r}$ by the martingale CLT. Lemmas \ref{ineq} and \ref{full-ineq} are used to simplify the expression of the martingale difference and verify Lindeberg's condition, respectively. The (limiting) covariance function is calculated based on Lemma \ref{double-e} with the help of Lemma \ref{full-ineq}. A new finding here is that the nonlinear dependency comes up with an extra term in the covariance function (Lemma \ref{double-e}), which results in a novel procedure of proving the convergence of this term.
	\item[] {\bf Step 3}: Tightness of $M_n^{(1)}(z)$ on $\mathcal C_n$. This step presents the basic idea for establishing the tightness. A key element is the uniform boundedness  of $E||(B_n-zI)^{-q}||$ for $q>0$ which is derived by Lemma \ref{lambda-bound}. By virtue of this and Lemmas \ref{ineq}-\ref{lambda-bound}, the tightness can be verified following similar arguments in \cite{BS04}.
	\item[] {\bf Step 4}: Convergence of $M_n^{(2)}(z)$. This final step mainly calculates the limit of $M_n^{(2)}(z)$, the limiting mean function of the LSSs. Akin to deriving the covariance function in Step 2, the nonlinear effect again brings an additional quantity to the mean function. Hence, most work in this part is to handle this new quantity. Note that Lemma \ref{lambda-bound} is useful in this step to obtain several convergence results and uniform boundedness on $\mathcal{C}_n$.
\end{itemize}
These detailed  four steps are presented in the next subsection. We note that when simplifying $M_n^{(1)}(z)$ and $M_n^{(2)}(z)$, one more straightforward method is used (see \eqref{undMn} and \eqref{Qnz} respectively), which are different from \cite{BS04}.

\subsection{Truncation and rescaling of the $\xi$-variable}

From the moment condition $E|(\xi_1^2-p)/\sqrt{p}|^{2+\varepsilon}<\infty$ for some $\varepsilon>0$ in Assumption (b), we can choose a sequence of $\delta_n>0$ such that
\begin{equation}\label{delta-n}
\delta_n\rightarrow0,\quad \delta_nn^{1/2}\rightarrow\infty, \quad \delta_n^{-2}p^{-1}E[(\xi_1^2-p)^2I_{\{|\xi_1^2-p|\geq\delta_np\}}]\rightarrow0.
\end{equation}
Let $\hat B_n=\sum_{j=1}^n\hat \x_j\hat\x_j'/n$ where $\hat\x_j=\hat\xi_j A\bu_j$ with $\hat\xi_j=\xi_jI_{\{|\xi_j^2-p|<\delta_np\}}$.
We then have
\begin{eqnarray}
P(\hat B_n\neq B_n)&\leq& nP(|\xi_1^2-p|\geq\delta_np)\nonumber\\
&\leq& \delta_n^{-2}np^{-2}E[(\xi_1^2-p)^2I_{\{|\xi_1^2-p|\geq\delta_np\}}]\rightarrow0.\label{bn-hat}
\end{eqnarray}
Define $\tilde B_n=\sum_{j=1}^n\tilde\x_j\tilde\x_j'/n$ where $\tilde\x_j=\tilde\xi_jA\bu_j$ with $\tilde\xi_j=\hat\xi_j/\sigma_n$ and $\sigma_n^2=E(\hat\xi_1^2)/p$.
By the assumptions in \eqref{delta-n},
\begin{eqnarray}
|1-\sigma_n^2|&=&E(\xi_1^2/p-1)I_{\{|\xi_1^2-p|\geq\delta_np\}}+EI_{\{|\xi_1^2-p|\geq\delta_np\}}\nonumber\\
&\leq&\delta_n^{-1}(1+\delta_n^{-1})p^{-2}E(\xi_1^2-p)^2I_{\{|\xi_1^2-p|\geq\delta_np\}}=o(p^{-1}).\label{sigma-n}
\end{eqnarray}
Therefore we have
\begin{equation*}
E(\tilde \xi_1^2)=p\quad \text{and}\quad E(\tilde \xi_1^4)=\frac{1}{\sigma_n^4}\left(E(\xi_1^4)-E\xi_1^4I_{\{|\xi_1^2-p|\geq\delta_np\}}\right)=p^2+\tau p+o(p).
\end{equation*}
On the other hand, write $\bu_j=\z_j/||\z_j||$ where,  and in the following $\z_j\sim N(0,I_p)$  and $\|\cdot\|$ denotes the spectral norm for a matrix, or $L_2$ norm for a vector. By the strong law of large numbers, for any fixed $0<\eta<1$,  we have $\max\{||\z_j||^2/p: j=1,\ldots,n\}\geq 1-\eta$ holds almost surely for large $p$. Hence we have for large $p$,
$$
||\tilde B_n||=\bigg|\bigg|\frac{1}{n}\sum_{j=1}^n\frac{\tilde\xi_j^2/p}{||\z_j||^2/p}A\z_j\z_j'A'\bigg|\bigg|\leq \frac{(1+\delta_n)}{(1-\eta)\sigma_n^2}||\Sigma||\bigg|\bigg|\frac{1}{n}\sum_{j=1}^n\z_j\z_j'\bigg|\bigg|
$$
almost surely. Thus, from \cite{Yin88}, we know that $\limsup_n\lambda_{\max}^{\tilde B_n}$ (and similarly  $\limsup_n\lambda_{\max}^{\hat B_n}$) are almost surely bounded by $\lim\sup_p ||\Sigma||(1+\sqrt{c})^2$.

Let $\tilde G_n(x)$ and $\hat G_n(x)$ be the analogues of $G_n(x)$ with the matrix
$B_n$ replaced by $\tilde B_n$ and $\hat B_n$, respectively. From the arguments in \cite{BS04} and \eqref{sigma-n}, we can get for $f(x,z)=1/(x-z)$ $(z\in \mathcal C)$, almost surely,
\begin{align*}
	&p^2\bigg|\int f(x)d\hat G_n(x)-\int f(x)d\tilde G_n(x)\bigg|^2\leq \left(\sum_{j=1}^pK|\lambda_j^{\hat B_n}-\lambda_j^{\tilde B_n}|\right)^2\\
	&\leq4K^2\sum_{j=1}^p\left(\sqrt{\lambda_j^{\hat B_n}}-\sqrt{\lambda_j^{\tilde B_n}}\right)^2\sum_{j=1}^p\left(\lambda_j^{\hat B_n}+\lambda_j^{\tilde B_n}\right)\\
	&\leq4K^2p(\lambda_{\max}^{\hat B_n}+\lambda_{\max}^{\tilde B_k})\frac{1}{n}\sum_{j=1}^n(\hat\x_j-\tilde\x_j)'(\hat\x_j-\tilde\x_j)\\
	&=4K^2p(\lambda_{\max}^{\hat B_n}+\lambda_{\max}^{\tilde B_k})\frac{(1-\sigma_n^2)^2}{\sigma_n^2(1+\sigma_n^2)}\tr(\hat B_n)\rightarrow0,
\end{align*}
where $K$ is an upper bound of $|f'_x(x,z)|$. As a consequence of this and \eqref{bn-hat},
$$
M_n(z)=p\int f(x)d G_n(x)=p\int f(x)d\tilde G_n(x)+o_p(1).
$$
Therefore, we only need to find the limiting distribution of $\int f(x)d\tilde G_n(x).$
For simplicity, we still use $B_n$, $\x_j$, $\xi_j$ instead of $\tilde B_n$, $\tilde\x_j$, $\tilde \xi_j$, respectively, and assume that
\begin{equation}\label{xi-ass}
\forall j,\ |\xi_j^2-p|< \delta_np,\quad E\left(\xi_1^2\right)=p,\quad E\left(\xi_1^4\right)=p^2+\tau p+o(p),
\end{equation}
in the sequel.

\subsection{Finite dimensional convergence of $M_n^{(1)}(z)$ in distribution}
We will show in this part that for any positive integer $r$
and any complex numbers
$z_{1},\ldots, z_r \in \mathcal C_n$, the random vector
$$[M_{n}^{(1)}(z_1),\ldots, M_n^{(1)}(z_r)]$$
converges to a $2r$-dimensional Gaussian vector. Because of
Assumption (c), without loss of generality, we may assume $\|\Sigma\|\leq1$ for all $p$.
We will denote by $K$ any constant appearing in inequalities and it may take different values at different
places.

We first define some quantities which are frequently used in the sequel:
\begin{align*}
	&r_j=(1/\sqrt{n})\x_j,\quad D(z)=B_n-zI,\quad D_j(z)=D(z)-r_jr'_j,\\
	 &D_{ij}(z)=D(z)-r_ir_i'-r_jr'_j,\quad
	\varepsilon_j(z)=r_j'D_j^{-1}(z)r_j-\frac{1}{n}\tr\Sigma D_j^{-1}(z),\\  &\zeta_j(z)=r_j'D_j^{-2}(z)r_j-\frac{1}{n}\tr\Sigma D_j^{-2}(z),\quad
\beta_j(z)=\frac{1}{1+r_j'D_j^{-1}(z)r_j},\\  &\bar\beta_j(z)=\frac{1}{1+n^{-1}\tr\Sigma D_j^{-1}(z)},\quad b_n(z)=\frac{1}{1+n^{-1}E\tr\Sigma D_j^{-1}(z)}.
\end{align*}
Note that, for any $z=u+iv\in \mathbb C^+$, the last three quantities are bounded in absolute value by $|z|/v$. Moreover,  $D^{-1}(z)$ and $D_j^{-1}(z)$ satisfy
\begin{equation}\label{D-D}
D^{-1}(z)-D_j^{-1}(z)=-D_j^{-1}(z)r_jr_j'D_j^{-1}(z)\beta_j(z).
\end{equation}
From Lemma 2.6 in \cite{SB95}, for any $p\times p$ matrix $B$,
\begin{equation}\label{d-d}
|\tr(D^{-1}(z)-D_j^{-1}(z))B|\leq\frac{||B||}{v}.
\end{equation}

Let $E_0(\cdot)$ denote expectation, and $E_j(\cdot)$ the conditional given the $\sigma$-field generated by $r_1,\ldots,r_j$.
Using the martingale decomposition, we can express $M_n^{(1)}(z)$ as
\begin{align*}
&\sum_{j=1}^{n}(E_j-E_{j-1})\tr D^{-1}(z)
	=\sum_{j=1}^n (E_j- E_{j-1})\tr\left[D^{-1}(z)-D_j^{-1}(z)\right]\\
	=&-\sum_{j=1}^n(E_j-E_{j-1})\beta_j(z)r_j'D_j^{-2}r_j
	=\sum_{j=1}^n(E_j-E_{j-1})\frac{d\log(\beta_j(z)/\bar{\beta}_j(z))}{d z},
\end{align*}
where the second equality uses the identity \eqref{D-D}. By the fact that  $$\beta_j(z)=\bar\beta_j(z)-\bar\beta_j(z)\beta_j(z)\varepsilon_j(z)=\bar \beta_j(z)-\bar\beta_j^2\varepsilon_j(z)+\bar\beta_j^2(z)\beta_j(z)\varepsilon_j^2(z),$$ we have
\begin{align}\label{undMn}
	M_n^{(1)}(z)=\frac{d}{d z}\sum_{j=1}^n(E_j-E_{j-1})\log[1-\bar\beta_j(z)\varepsilon_j(z)+\bar\beta_j(z)\beta_j(z)\varepsilon_j^2(z)].
\end{align}
Notice that for all $j>0$ and  any $z\in{\cal C}_n$, $\bar\beta_j(z)\varepsilon_j(z)$  and $\bar\beta_j(z)\beta_j(z)\varepsilon_j^2(z)$ are almost surely away from 1 when $n$ is large enough. In addition,
by  Lemma  \ref{ineq} and  Burkholder's inequality (Lemma  2.1 in \cite{BS04}), we   have
\begin{eqnarray*}
	E\bigg|\sum_{j=1}^n(E_j-E_{j-1})\bar\beta_j(z)\beta_j(z)\varepsilon^2_j(z)\bigg|^2
	=O(\de_n^2)\to0.
\end{eqnarray*}
Therefore, applying  Taylor expansion,
\begin{align*}
M_n^{(1)}(z)=&-\frac{d}{dz}\sum_{j=1}^n (E_j-E_{j-1})\bar\beta_j(z)\varepsilon_j(z)+o_p(1)\\
=&-\frac{d}{dz}\sum_{j=1}^n E_j(\bar\beta_j(z)\varepsilon_j(z))+o_p(1).
\end{align*}
For any $\epsilon>0$,
\begin{eqnarray*}
&&	\sum_{j=1}^nE\left|E_j\frac{d}{dz}\varepsilon_{j}(z)\bar\beta_{j}(z)\right|^2 I_{\left(\left|E_j\frac{d}{dz}\varepsilon_{j}(z)\bar\beta_{j}(z)\right|\geq\epsilon\right)}
	\leq\frac{1}{\epsilon^2}\sum_{j=1}^n E\Big|E_j\frac{d}{dz}\varepsilon_{j}(z)\bar\beta_{j}(z)\Big|^4\\
	&\leq&\frac{K}{\varepsilon^2}\sum_{j=1}^n
	\left(\frac{|z|^4E|\zeta_j(z)|^4}{v^4}+\frac{|z|^8p^4E|\varepsilon_j(z)|^4}{v^{16}n^4}\right)
\end{eqnarray*}
which tends to zero according to Lemma \ref{full-ineq}, and thus Lindeberg's condition is verified.
Therefore, from the martingale CLT (Theorem 35.12 \cite{B95}),
the random vector $(M_{n}^{(1)}(z_{j}))$ tends to a $2r$-dimensional zero-mean Gaussian vector
$(M(z_{j}))$ with covariance function $Cov(M(z_1),M(z_2))$ being
\begin{equation}\label{lcov}
\lim_{n\to\infty}\sum_{j=1}^n \frac{\partial^2}{\partial z_1\partial z_2}E_{j-1}\left(E_j\varepsilon_{j}(z_1)\bar\beta_{j}(z_1)\cdot
E_j\varepsilon_{j}(z_2)\bar\beta_{j}(z_2)\right),
\end{equation}
provided that this limit exits in probability. The same argument in page 571 of \cite{BS04} implies that it suffices to show
\begin{equation}\label{lcov-1}
\sum_{j=1}^nE_{j-1}\prod_{k=1}^2E_j\bar \beta_j(z_k)\varepsilon_j(z_k)
\end{equation}
converges in probability. In addition, by the martingale decomposition,
\begin{align}
E|\bar\beta_j(z)-b_n(z)|^2
=&|b_n(z)|^2n^{-2}E\bigg|\bar\beta_1(z)\sum_{k=2}^n(E_k-E_{k-1})\tr(D_1^{-1}(z)-D_{1k}^{-1}(z))\bigg|^2\nonumber\\
\leq&K|z|^4v^{-6}n^{-1},\label{beta-b}
\end{align}
where the inequality is from \eqref{d-d}.
 Thus it is sufficient to
study the convergence of
\begin{equation}\label{cov-term-2}
b_n(z_1)b_n(z_2)\sum_{j=1}^nE_{j-1}\left(E_j\varepsilon_j(z_1)E_j\varepsilon_j(z_2)\right),
\end{equation}
whose second mixed partial derivative yields the limit of \eqref{lcov}.
Applying Lemma \ref{double-e}, we know that
\begin{equation}\label{cov-term-3}
\eqref{cov-term-2}= n\left(\frac{E\xi^4}{p(p+2)}-1\right)T_1+\frac{2E\xi^4}{p(p+2)}T_2,
\end{equation}
where
\begin{eqnarray*}
	&&T_1=b_n(z_1)b_n(z_2)\frac{1}{n^3}\sum_{j=1}^n\tr[\Sigma E_jD_j^{-1}(z_1)]\tr[\Sigma E_jD_j^{-1}(z_2)],\\
	&&T_2=b_n(z_1)b_n(z_2)\frac{1}{n^2}\sum_{j=1}^n\tr[\Sigma E_jD_j^{-1}(z_1)\Sigma E_jD_j^{-1}(z_2)].
\end{eqnarray*}

We note that the statistic $T_2$ has the same form as Equation (2.8) in \cite{BS04}, which has been handled under their model.
Following their calculations and using Lemmas \ref{ineq}-\ref{full-ineq} instead, one may get
\begin{eqnarray}\label{limT2}
	T_2
	\xrightarrow{i.p.}\int_0^{a(z_1,z_2)}\frac{1}{1-z}dz,
\end{eqnarray}
where
\begin{eqnarray*}
	a(z_1,z_2)=\int\frac{c\um(z_1)\um(z_2)t^2dH(t)}{(1+t\um(z_1))(1+t\um(z_2))}=1+\frac{\um(z_1)\um(z_2)(z_1-z_2)}{\um(z_2)-\um(z_1)},
\end{eqnarray*}
and
\begin{equation}\label{p-t2}
\frac{\partial^2T_2}{\partial z_1\partial z_2}	\xrightarrow{i.p.}\frac{\um'(z_1)\um'(z_2)}{(\um(z_1)-\um(z_2))^2}-\frac{1}{(z_1-z_2)^2}.
\end{equation}

Now we derive the limit of $T_1$ and its second mixed partial derivative, which is new compared with the linear transform model.  Denote
\begin{eqnarray*}
	\beta_{ij}(z)={(1+r_i'D_{ij}^{-1}(z)r_i)^{-1}},\quad b_{1}(z)=({1+n^{-1}E\tr \Sigma D_{12}^{-1}(z)})^{-1}.
\end{eqnarray*}
By similar proofs of  \eqref{beta-b}  and Equation (4.3) of \cite{BS98}, one may get
$
|b_1(z)-b_n(z)|\leq Kn^{-1}$ and $ |b_n(z)-E\beta_1(z)|\leq Kn^{-1/2},
$
respectively. Also, by Equation  (2.2) of \cite{S95} and discussions in Section 5 of \cite{BS98}, we obtain
$$
E\beta_1(z)=-zE\underline{m}_n(z)\quad\text{and}\quad |E\um_n(z)-\um_{F^{c_n,H_p}}(z)|\leq Kn^{-1},
$$
respectively. Therefore, we get
\begin{equation}\label{b-zm}
|b_1(z)+z\um_{F^{c_n,H_p}}(z)|\leq Kn^{-1/2}.
\end{equation}
With the quantity $b_1(z)$, we define a nonrandom matrix $L(z)$ for the purpose of replacing $D_j(z)$ in $T_1$,
\begin{align*}
L(z)=-zI+\frac{n-1}{n}b_1(z)\Sigma,
\end{align*}
which satisfies
\begin{equation}\label{z-b}
||L(z)||^{-1}\leq \frac{|b_1^{-1}(z)|}{\Im(zb_1^{-1}(z))}\leq\frac{|b_1^{-1}(z)|}{\Im(z)}\leq\frac{1+p/(nv)}{v}.
\end{equation}
By the identity $r_i'D_j^{-1}(z)=\beta_{ij}(z)r_i'D_{ij}^{-1}(z)$, we get their difference
\begin{eqnarray}
D_j^{-1}(z)-L^{-1}(z)=b_1(z)R_1(z)+R_2(z)+R_3(z),\label{dj-1}
\end{eqnarray}
where
\begin{align*}
	R_1(z)=&-\sum_{i\neq j}L^{-1}(z)(r_ir_i'-n^{-1}\Sigma)D_{ij}^{-1}(z),\\
	R_2(z)=&-\sum_{i\neq j}(\beta_{ij}(z)-b_1(z))L^{-1}(z)r_ir_i'D_{ij}^{-1}(z),\\
	R_3(z)=&-n^{-1}b_1(z)L^{-1}(z)\Sigma \sum_{i\neq j}\left(D_{ij}^{-1}(z)-D_{j}^{-1}(z)\right).
\end{align*}
For any $p\times p$ (nonrandom) matrix $M$, 
from \eqref{d-d}, \eqref{z-b}, and Lemma \ref{full-ineq}, we get
\begin{align}
E|\tr R_1(z)M|\leq& nE^{1/2}|r_1'D_{12}^{-1}(z)ML^{-1}(z)r_1-n^{-1}\tr\Sigma D_{12}^{-1}(z)ML^{-1}(z)|^2\nonumber\\
\leq&n^{1/2}K||M||\frac{(1+p/(nv))}{v^2}\label{AM},
\end{align}
\begin{align}
E|\tr R_2(z)M|\leq& nE^{1/2}(|\beta_{12}(z)-b_1(z)|^2)E^{1/2}\bigg|r_1'D_{12}^{-1}ML^{-1}(z)r_1\bigg|^2\nonumber\\
\leq&n^{1/2} K ||M||\frac{|z|^2(1+p/(nv))}{v^5},\label{BM}\\
|\tr R_3(z)M|&\leq||M||\frac{|z|(1+p/(nv))}{v^3}.\label{CM}
\end{align}
Hence, plugging \eqref{dj-1} into $T_1$ and applying the inequalities \eqref{b-zm},  \eqref{AM}-\eqref{CM}, we have
\begin{align*}
	\prod_{k=1}^2\tr E_jD_j^{-1}(z_k)\Sigma=&\prod_{k=1}^2\tr L^{-1}(z_k)\Sigma+Q_1(z_1,z_2)\\
	=&p^2\prod_{k=1}^2\frac{1}{z_k}\int\frac{tdH_p(t)}{1+t\um_{F^{c_n,H_p}}(z_k)}+Q_2(z_1,z_2),
\end{align*}
where $E|Q_k(z_1,z_2)|\leq Kn^{3/2}$, $k=1,2$. We thus get
$$
T_1=\prod_{k=1}^2\um_{F^{c_n,H_p}}(z_k)\int\frac{c_ntdH_p(t)}{1+t\um_{F^{c_n,H_p}}(z_k)}+o_p(1)\xrightarrow{i.p.}\prod_{k=1}^2\left(1+z_k\um(z_k)\right)
$$
whose second mixed partial derivative is
\begin{equation}\label{p-t1}
{\partial^2T_1}/({\partial z_1\partial z_2})\xrightarrow{i.p.}(\um(z_1)+z_1\um'(z_1))(\um(z_2)+z_2\um'(z_2)).
\end{equation}
The result in \eqref{p-t1} can be obtained by Vitali's convergence theorem (Lemma 2.3 in \cite{BS04}).

Collecting results in \eqref{cov-term-3}, \eqref{p-t2} and \eqref{p-t1}, we finally get
\begin{align*}
	Cov(M(z_1),M(z_2))=&(\um(z_1)+z_1\um'(z_1))(\um(z_2)+z_2\um'(z_2))\\
	&+{2\um'(z_1)\um'(z_2)}/{(\um(z_1)-\um(z_2))^2}-{2}/{(z_1-z_2)^2},
\end{align*}
which completes  the  proof
 of Step 1.
\subsection{Tightness of $M_n^{(1)}(z)$}

The tightness of $M_n^{(1)}(z)$ can be established by verifying the moment condition (12.51) of \cite{B68}, i.e.,
\begin{equation}\label{tightness}
\sup_{n,z_1,z_2\in \mathcal C_n}{E|M_n^{(1)}(z_1)-M_n^{(1)}(z_2)|^2}/{|z_1-z_2|^2}<\infty.
\end{equation}
By the martingale decomposition and  the equality
\begin{align*}
m_n(z_1)-m_n(z_2)={(z_1-z_2)}{p^{-1}}\tr(D^{-1}(z_1)D^{-1}(z_2)),
\end{align*}
to show \eqref{tightness}, it is sufficient to prove  the absolute second moment of
\begin{align*}
&\sum_{j=1}^n(E_j-E_{j-1})\tr[D^{-1}(z_1)D^{-1}(z_2)]
\end{align*}
is bounded uniformly. We first show the uniformly boundedness of $E||D^{-q}(z)||$ on $\mathcal C$ for any fixed $q>0$. Note that $D^{-1}(z)$ is bounded on $z\in \mathcal C_u$. While for $z\in \mathcal C_l\cup C_r$, applying Lemma \ref{lambda-bound} with suitable large $s$, we have uniformly
\begin{align*}
	E||D^{-1}(z)||^q\leq K+\frac{1}{v^q}P(||B_n||>\eta_r\ \text{or}\ \lambda_{\min}^{B_n}<\eta_l)
	\leq K+o(1),
\end{align*}
where $\lim\sup_p ||\Sigma||(1+\sqrt{c})^2<\eta_r<x_r$ and $x_l<\eta_l<\lim \inf_p \lambda_{\min}^{\Sigma}(1-\sqrt{c})^2$. Analogously $E||D_j^{-1}(z)||^q$ has the same order, and we thus get
\begin{equation}\label{D-bound}
\max \{ E||D^{-1}(z)||^q, E||D^{-1}_j(z)||^q, E||D^{-1}_{ij}(z)||^q\}\leq K_q.
\end{equation}
Then, \eqref{tightness} can be obtained by the same procedure in Section 3 of  \cite{BS04}, applying Lemmas \ref{ineq}-\ref{lambda-bound} together with \eqref{D-bound}. We omit the details.

\subsection{Convergence of $M_n^{(2)}(z)$}
Next we will show that for $z \in\mathcal C_n$, $\{M_n^{(2)}(z)\}$  converges  to \eqref{mean},
is bounded and   forms a uniformly equicontinuous family.

We first introduce some auxiliary results, which can be verified by applying Lemma \ref{lambda-bound} in our theoretical framework through a similar proof of the same statements in \cite{BS04}.  First of all, we note that
\begin{eqnarray}\label{supsup-bound}
\sup_{z\in \mathcal C_n}|E\um_n(z)-\um(z)|\rightarrow0\quad \text{and}\quad \sup_{n, z\in \mathcal C_n}||V^{-1}(z)||< \infty,
\end{eqnarray}
where $V(z)=E\um_n(z)\Sigma+I.$
Then, for any nonrandom $p\times p$ matrix $M$,
\begin{eqnarray}\label{tr-dm-bound}
E|\tr D^{-1}(z)M-E\tr D^{-1}(z)M|^2\leq K||M||^2.
\end{eqnarray}
Next, there exists a number $\theta\in (0,1)$ such that for all $n$ large enough
\begin{eqnarray}\label{int-bound}
\sup_{z\in \mathcal C_n}\bigg|c_n\int \frac{\left(tE\um_n(z)\right)^2}{(1+tE\um_n(z))^2}dH_p(t)\bigg|<\theta.
\end{eqnarray}
Lastly, from (4.12) of \cite{BS04} and (5.2) in \cite{BS98}, we have  that
\begin{equation}\label{em-m0}
M_n^{(2)}(z)=-\frac{\um_{F^{c_n,H_p}}(z)Q_n(z)}{\left(1-\int\frac{c_nE\um_n(z)\um_{F^{c_n,H_p}}(z)t^2dH_p(t)}{(1+tE\um_n(z))(1+t\um_{F^{c_n,H_p}}(z))}\right)},
\end{equation}
where
\begin{align}\label{Qnz}
	Q_n(z)=&n\left(c_n\int\frac{dH_p(t)}{1+tE\um_n(z)}+zc_nEm_n(z)\right)\nonumber\\
	=&nE\beta_1(z)\left(r_1'D_1^{-1}(z)V^{-1}(z)r_1-n^{-1}E\tr V^{-1}(z)\Sigma D^{-1}(z)\right).
\end{align}
From \eqref{int-bound} and an analog inequality involving $\um_{F^{c_n,H_p}}(z)$, the denominator of \eqref{em-m0} is bounded away from zero.
Therefore, we need only to study the limit of $Q_n(z)$ for $z\in \mathcal C_n$.

For simplicity, we suppress the variable $z$ from expressions in the sequel when there is no confusion.
Let $\varrho_1:=\varrho_1(z)=r_1'D_1^{-1}r_1-(1/n)E\tr\Sigma D_1^{-1}.$
By the equality  \begin{align}\label{betab}
\beta_1=b_n-b_n\beta_1\varrho_1=b_n-b^2_n\varrho_1+b^2_n\beta_1\varrho^2_1,
\end{align}we have
$
Q_n=Q_n^{(1)}+Q_n^{(2)}+Q_n^{(3)},
$
where
\begin{align*}
Q_n^{(1)}&=b_nE(\tr D_1^{-1}V^{-1}\Sigma-\tr V^{-1}\Sigma D^{-1}),\\
Q_n^{(2)}&=-nb_n^2E\varrho_1(r_1'D_1^{-1}V^{-1}r_1-n^{-1}\tr D_1^{-1}V^{-1}\Sigma),\\
Q_n^{(3)}&=nb_n^2E\beta_1\varrho_1^2\left(r_1'D_1^{-1}V^{-1}r_1-n^{-1}E\tr V^{-1}\Sigma D^{-1}\right).
\end{align*}
For  $Q_n^{(1)}$, apply  \eqref{D-D} and \eqref{betab} again,
\begin{align}
	&E\tr V^{-1}\Sigma (D_1^{-1}-D^{-1})
	=E\beta_1r_1'D_1^{-1} V^{-1}\Sigma D_1^{-1}r_1\nonumber\\
	=&b_nn^{-1}E\tr D_1^{-1}V^{-1}\Sigma D_1^{-1}-b_nE\beta_1\varrho_1r_1'D_1^{-1} V^{-1}\Sigma D_1^{-1}r_1.\label{e-dd}
\end{align}
By Lemma \eqref{ineq}, H{\"o}lder's inequality and the fact that $r_1'D_1^{-1} V^{-1}\Sigma D_1^{-1}r_1$, $b_n$ and $\beta_1$   are all bounded for $z\in \mathcal C_n$, the second term in Equation \eqref{e-dd} is $o(1)$. Analogously,  we  can get that $Q_n^{(3)}=o(1)$.
Together with applying Lemma \ref{double-e} to $Q_n^{(1)}$, we finally obtain that
\begin{align}
Q_n
=&-b_n^2n^{-1}\left(\frac{E\xi^4}{p(p+2)}-1\right)E\tr D_1^{-1}\Sigma E\tr D_1^{-1}V^{-1}\Sigma\nonumber\\
&-b_n^2n^{-1}\left(\frac{2E\xi^4}{p(p+2)}-1\right)E\tr D_1^{-1}V^{-1}\Sigma D_1^{-1}\Sigma+o(1)\nonumber\\
:=&-(\tau-2)c_n^{-1}b_n^2Q_n^{(4)}-b_n^2Q_n^{(5)}+o(1),\label{t34}
\end{align}
where
$Q_n^{(4)}=n^{-2}E\tr D^{-1}\Sigma E\tr D^{-1}V^{-1}\Sigma$ and $Q_n^{(5)}=n^{-1}E\tr D^{-1}V^{-1}\Sigma D^{-1}\Sigma$.
The limit of $Q_n^{(5)}$ can be obtained by a similar approach to deriving (4.13)- (4.22) in \cite{BS04}. It turns out that
\begin{align}
Q_n^{(5)}=\frac{c_n}{z^2}\int\frac{t^2dH_p(t)}{(1+tE\um_n)^3}\left( 1-c_n\int\frac{(tE\um_n)^2dH_p(t)}{(1+tE\um_n)^2}\right)^{-1}+o(1).\label{t4}
\end{align}

The quantity $Q_n^{(4)}$ is new under the elliptical model. To study its limit, similar to \eqref{dj-1}, we approximate the matrix $D^{-1}(z)$ by
$$
\tL=-zI+b_n\Sigma.
$$
Notice that
$$
b_n=E\beta_1+O(n^{-1/2})=-zE\um_n+O(n^{-1/2})\rightarrow -z\um,
$$
as $n\rightarrow\infty$. By  \eqref{supsup-bound},
it follows that $\tL$ is nonsingular and $||\tL^{-1}||$ is bounded.
Then,  analogous to \eqref{dj-1}-\eqref{CM}, we have
\begin{eqnarray}
D^{-1}-\tL^{-1}
=b_n\R_1+\R_2+\R_3,\label{tilde-D-L}
\end{eqnarray}
where
\begin{align}
&|E\tr \R_1M|\leq n^{1/2}K,\quad \label{tilde-AM}
|E\tr \R_2M|\leq n^{1/2} K (E||M||^4)^{1/4},\\
&|\tr \R_3M|\leq K(E||M||^2)^{1/2}\label{tilde-CM}
\end{align}
for any $p\times p$ nonrandom matrix $M$ with bounded norm.
From \eqref{tilde-D-L}-\eqref{tilde-AM}, we have that
\begin{align}
n^{-1}E\tr D^{-1}\Sigma
&=-\frac{c_n}{z}\int\frac{tdH_p(t)}{1+tE\um_n}+o(1),\label{mean-p1}\\
n^{-1}E\tr D^{-1}V^{-1}\Sigma
&=-\frac{c_n}{z}\int\frac{tdH_p(t)}{(1+tE\um_n)^2}+o(1),\label{mean-p2}
\end{align}
Equations \eqref{mean-p1} and \eqref{mean-p2} imply that
\begin{equation}\label{t3}
Q_n^{(4)}=\frac{c_n^2}{z^2}\int\frac{tdH_p(t)}{1+tE\um_n}\int\frac{tdH_p(t)}{(1+tE\um_n)^2}+o(1).
\end{equation}
Combining \eqref{em-m0}, \eqref{t34}, \eqref{t4} and \eqref{t3}, we finally get
\begin{align*}
	M_n^{(2)}(z)=&\bigg[(\tau-2)c_n\int\frac{t\um_{F^{c_n,H_p}}dH_p(t)}{1+tE\um_n}\int\frac{t\um_n^2dH_p(t)}{(1+tE\um_n)^2}\\
	&+c_n\int\frac{t^2\um_{F^{c_n,H_p}}\um_n^2dH_p(t)}{(1+tE\um_n)^3}\left( 1-c_n\int\frac{(tE\um_n)^2dH_p(t)}{(1+tE\um_n)^2}\right)^{-1}\bigg]\\
	&\times\left(1-c_n\int\frac{E\um_n\um_{F^{c_n,H_p}}t^2dH_p(t)}{(1+tE\um_n)(1+t\um_{F^{c_n,H_p}})}\right)^{-1}+o(1)\\
	\rightarrow&(\tau-2)\int\frac{(z\um+1)\um'tdH(t)}{(1+t\um)^2}+c\int\frac{(\um't)^2dH(t)}{\um(1+t\um)^3},
\end{align*}
where $\um'=\um'(z)$ denotes the derivative of $\um(z)$ with respective to $z$.

The boundedness and uniform equicontinuity  for $z\in \mathcal C_n$ can be verified directly following the arguments on Pages 592-593 of \cite{BS04}. So we omit them here. Then, the proof of Theorem \ref{clt} is complete.

\appendix

\section{Appendix}

These lemmas can be viewed as extensions of independent cases. Their proofs are postponed to  the supplementary material \citep{suppm}.

\begin{lemma}\label{double-e}
	Let $\x=\xi\bu$ where $\xi$ and $\bu$ are defined in Assumption (b). Then, for any $p\times p$ complex matrices $C$ and $\tilde C$,
	\begin{equation}
	E(\x'C\x-\tr C)(\x'\tilde C\x-\tr\tilde C)=\frac{E\xi^4}{p(p+2)}\left(\tr C\tr\tilde C+\tr C\tilde C'+\tr C\tilde C\right)-\tr C\tr\tilde C.
	\end{equation}
\end{lemma}
\begin{lemma}\label{ineq}
	Let $\x=\xi\bu$ where $\xi$ satisfies \eqref{xi-ass}, independent of $\bu\sim U(S^{p-1})$, then for any $p\times p$ complex matrix $C$ and  $q\geq 2$,
	\begin{equation}\label{q-moment}
	E|\x'C\x-\tr C|^q\leq K||C||^q\delta_n^{q-2}p^{q-1},
	\end{equation}
	where $K$ is a positive constant depending only on $q$.
\end{lemma}

\begin{lemma}\label{full-ineq}
	Let $r=\xi \bu/\sqrt{n}$ where $\xi$ satisfies \eqref{xi-ass}, independent of $\bu\sim U(S^{p-1})$. Then, for any nonrandom $p\times p$ matrix $C_k$, $k=1,\ldots,q_1$ and $\tilde C_l$, $l=1,\ldots,q_2$, $q_1, q_2\geq0$,
	\begin{equation*}
	\bigg|E\left(\prod_{k=1}^{q_1}r'C_kr\prod_{l=1}^{q_2}(r'\tilde C_lr-n^{-1}\tr\tilde C_l)\right)\bigg|
	\leq Kn^{-(1\wedge q_2)}\delta_n^{(q_2-2)\vee0}\prod_{k=1}^{q_1}||C_k||\prod_{l=1}^{q_2}||\tilde C_l||,
	\end{equation*}
	where $K$ is a positive constant depending on $q_1$ and $q_2$.
\end{lemma}

\begin{lemma}\label{lambda-bound}
	Suppose \eqref{xi-ass} holds. Then, for any positive $s$,
	$$
	P(||B_n||>\eta_r)=o(n^{-s}),
	$$
	whenever $\eta_r>\lim\sup_{p\rightarrow\infty} ||\Sigma||(1+\sqrt{c})^2$.
	If $0<\lim\inf_{p\rightarrow\infty}\lambda_{\min}^{\Sigma}I_{(0,1]}(c)$ then,
	$$
	P(\lambda_{\min}^{B_n}<\eta_l)=o(n^{-s}),
	$$
	whenever $0<\eta_l<\lim\inf_{p\rightarrow\infty}\lambda_{\min}^{\Sigma}I_{(0,1)}(c)(1-\sqrt{c})^2$.
\end{lemma}

\section*{Acknowledgement}
The authors are grateful to Prof.\ Zhidong Bai for discussions on the innovation of high-dimensional elliptical models, and to Prof.\ Jianfeng Yao for discussions on spiked covariance models. They also would like to thank the editors and referees whose careful and detailed comments have led to great improvements of the manuscript.


\end{document}